\documentclass[11pt,a4paper]{article}

\usepackage{epsf,epsfig,amsfonts,amsgen,amsmath,amstext,amsbsy,amsopn,amsthm
}
\usepackage{amsmath,times,mathptmx}
\usepackage{amsfonts,amsthm,amssymb}
\usepackage{amsfonts}
\usepackage{graphics}
\usepackage{latexsym,bm}
\usepackage{amsfonts,amsthm,amssymb,bbding}
\usepackage{indentfirst}
\usepackage{graphicx}
\usepackage{color}
\usepackage[colorlinks=true,anchorcolor=blue,filecolor=blue,linkcolor=blue,urlcolor=blue,citecolor=blue]{hyperref}
\usepackage{float}
\usepackage{tikz}
\newtheorem{thm}{Theorem}

\newtheorem{lemma}{Lemma}
\newtheorem{false statement}{False statement}

\newtheorem{fact}{Fact}

\theoremstyle{definition}

\newtheorem{claim}{Claim}

\newtheorem{conj}{Conjecture}

\newtheorem{prob}{Problem}[section]

\begin{document}

\title{Spectral radius, edge-disjoint cycles and cycles of the same length
\footnote{Supported by the National Natural Science Foundation of China (11771141, 12171066 and 12011530064).}}
\author{{Huiqiu Lin$^a$, Mingqing Zhai$^b$ and Yanhua Zhao$^a$}\thanks{E-mail: huiqiulin@126.com (H.Q. Lin), mqzhai@chzu.edu.cn (M. Zhai), yhua030@163.com (Y.H. Zhao)}\\
{\footnotesize $^a$School of Mathematics, East China University of Science and Technology, Shanghai 200237, China}\\
{\footnotesize $^b$ School of Mathematics and Finance, Chuzhou University, Chuzhou, Anhui 239012, China}}
\date{}

\maketitle {\flushleft\large\bf Abstract:}
In this paper, we give spectral conditions to guarantee the existence of two edge disjoint cycles and two cycles of the same length. These two results can be seen as spectral analogues of  Erd\H{o}s and Posa's size condition and Erd\H{o}s' classic problem on non existence of two cycles of the same length. By using double leading eigenvectors skill, we further give spectral condition to guarantee the existence of $k$ edge disjoint triangles.

\vspace{0.1cm}
\begin{flushleft}
\end{flushleft}

\section{Introduction}
The study of cycles has a rich history in graph theory. A
folklore result states that every
graph with minimum degree at least 2 has a cycle.
On the other hand, a connected graph on $n$ vertices without any cycle is a tree.
Let ${\Omega_k}'$ be the family of graphs with $k$
edge-disjoint cycles, and $\overline{\Omega'_k}$ be the family of graphs out of ${\Omega_k}'$.
Erd\H{o}s and Posa proved the following theorem.
\begin{thm}[Erd\H{o}s and Posa \cite{EP65}]\label{edcycle}
Let $G$ be a graph on $n$ vertices and $m$ edges. If $G\in \overline{\Omega'_2}$, then $m\leq n+3$.
Furthermore, if $G\in \overline{\Omega'_2}$ and $m=n+3$, then $G$ is obtained from a subdivision
$G_0$ of $K_{3,3}$ by adding a forest and exactly one edge, joining each
tree of the forest to $G_0$.
\end{thm}

Bollob\'{a}s in his classic book \cite{B78} posed a more general problem as follows:
What is the maximum size of  a graph $G\in \overline{\Omega'_k}$ of order $n$?
Up to now, this problem is still widely open.

Another classic problem involves cycles of the same length.
Let $f(n)$ be the maximum number of edges in a graph on $n$ vertices without two cycles of the same length.
Erd\H{o}s posed the problem of determining $f(n)$,
which was listed as one of $50$ unsolved problems in the textbook by Bondy and Murty
(see \cite[ p.247, Problem 11]{BM76}). It has attracted the attention of many scholars.
In 1988, Shi \cite{S88} gave a lower bound of $f(n)$, which was improved by Lai \cite{L17} in 2017.
Boros, Caro, F\"{u}redi and Yuster \cite{BC01} showed that $\sqrt{n}-o(\sqrt{n})\leq f(n)-n\leq 1.98\sqrt{n}(1+o(1))$,
and they further conjectured that $\lim\limits_{n\rightarrow\infty}\frac{f(n)-n}{\sqrt{n}}=1$.
Very recently, Ma and Yang \cite{MY} confirmed their conjecture by showing that $f(n)<n+\sqrt{n}+o(n)$ for any $n$-vertex 2-connected graph.
However, the exactly value of $f(n)$ is still unknown.
For more related results, we refer the interested readers to \cite{CL98,L03,L20}.

One main purpose of this paper is to investigate these two problems from a spectral perspective.
The eigenvalue conditions for cycles have been studied by a plenty of researchers (see \cite{ZL02})
and two surveys \cite{CZ18,N11} for more fruitful results.
In 1995, Favaron, Mah\'{e}o and Sacl\'{e} \cite{FMS93}  proved
that every graph on $n$ vertices with $\rho(G)>\rho(S_{n,1})=\sqrt{n-1}$
contains a $C_3$ or a $C_4$. Generalizing this theorem,
Nikiforov \cite{N10} conjectured that: (a)  every graph of sufficiently large order $n$ with $\rho(G)\geq \rho(S_{n,k})$
contains a $C_{2k+1}$ or a $C_{2k+2}$,
unless $G=S_{n,k}$.
He \cite{N10} also conjectured that: (b) every graph of sufficiently large order $n$ with $\rho(G)\geq \rho(S^+_{n,k})$
contains a $C_{2k+2}$, unless $G=S^+_{n,k}$.
The conjecture (a) was confirmed for $k=2$ in \cite{YWZ12} and conjecture (b)
was confirmed recently for $k=2$ in \cite{ZL20}. For
$C_4$, Nikiforov \cite{N09} and Zhai and Wang \cite{ZW12} characterized the extremal graphs on odd $n$ and even $n$, respectively.
For consecutive cycles, see Nikiforov \cite{N08}, Ning and Peng \cite{NP1}, Zhai and Lin \cite{ZL1} and Li and Ning \cite{LN}.
Only very recently, confirming the first open case of a conjecture due
to Bollob\'{a}s and Nikiforov, Lin, Ning and Wu \cite{LNW21+} obtained
a new eigenvalue condition for triangles. Furthermore, they obtained
a spectral analogue of a theorem of Erd\H{o}s, which states that a non-bipartite graph $G$ with $\rho(G)\geq \rho(S(K_{\lfloor\frac{n-1}{2}\rfloor,\lceil\frac{n-1}{2}\rceil}))$  contains a triangle unless $G\cong S(K_{\lfloor\frac{n-1}{2}\rfloor,\lceil\frac{n-1}{2}\rceil})$ where $S(K_{\lfloor\frac{n-1}{2}\rfloor,\lceil\frac{n-1}{2}\rceil})$
denotes a subdivision of $K_{\lfloor \frac{n-1}{2}\rfloor,\lceil \frac{n-1}{2}\rceil}$
on one edge.
There are also quite a lot of references on eigenvalues and large cycles, we refer the reader to \cite{GH19,GN20,LN}.

It seems that the study of cycles in terms of eigenvalues is closely related
to extremal graph theory. For example, in order to solve an open problem raised by Caro
and Yuster on degree powers in graphs with a forbidden even cycle, Nikiforov \cite{N09-EJC}
proved an extension of the classical Erd\H{o}s-Gallai
theorem on cycles. Interestingly, he also used it as a tool to study a spectral problem,
showing that $\rho^2(G)-k\rho(G)\leq k(n-1)$ if $G$ is $C_{2k}$-free.
Only very recently, it is also shown that the spectral theorem on
consecutive cycles \cite{N08} is a tool for studying some generalized
Tur\'an-type problems, such as estimating
ex$(n,P_{l},C_{2k+1})$ (see \cite{GGMV20}).

In this paper, we continue the project of studying cycles in spectral prospective.
Compared with \cite{LNW21+}, we focus
on eigenvalues conditions for edge-disjoint cycles, cycles
of the same length.

As usual, we denote by $K_r$ the complete graph on $r$ vertices,
by $S_{n, k}$ the graph obtained by joining $n-k$ isolated vertices
with each vertex in $K_k$, and by $K_{a, b}$ the complete bipartite
graph with two parts of orders $a$ and $b$, respectively.
Let $H\in \{C_r, K_r, K_{r,r}\}$ where $C_r$ denotes the cycle of length $r$. We denote by
$H\bullet K_{a, b}$ the graph obtained by coalescing a vertex of $K_{a,b}$ belonging to the part of size $a$ 
with a vertex in $H$ and retaining the
connection of edges in $H$ and $K_{a, b}$. Let $K_{1,n-1}^+$
be the graph obtained by adding an edge within the independent set of $K_{1,n-1}$.

It is rather surprising that we find a spectral analogue of $f(n)$ and the extremal graph is also determined.

\begin{thm}\label{thm1.2}
Let $G$ be a graph of order $n\geq 26$.
If $\rho(G)\geq \rho(K^+_{1,n-1})$, then $G$ contains two cycles with the same length unless $G\cong K_{1,n-1}^+$.
\end{thm}

In the following, we give a spectral analogue of Theorem \ref{edcycle}.
\begin{thm}\label{Thm:edgedisjointcycle}
Let $G$ be a graph of order $n\geq 17$.
If $\rho(G)\geq \rho(K_4\bullet K_{1,n-4})$,
then $G$ contains two edge-disjoint cycles unless $G\cong K_4\bullet K_{1,n-4}$.
\end{thm}
It is a natural wish to  determine the maximum spectral radius of $\Omega'_k$-free graphs of order $n$,
which is stated as follows.
\begin{prob}\label{prob0}
What is the maximum spectral radius of  $\Omega'_k$-free graphs of order $n$?
\end{prob}
However, it seems difficult to solve Problem \ref{prob0}, even for giving a conjecture on $k\geq 3$. Therefore, we turn out this problem in a special version, i.e., the Brualdi-Solheid-Tur\'{a}n type problem for $k$ edge-disjoint triangles graphs, which also has an interesting history. Let $\mathcal{F}$ be a family of graphs. A graph $G$ is called  $\mathcal{F}$-free if it contains no any graph in $\mathcal{F}$ as a subgraph. The Tur\'{a}n number of $\mathcal{F}$, denote by ex$(n,\mathcal{F})$, is the maximum number of edges in an $\mathcal{F}$-free graph of order $n$. Let EX$(n,\mathcal{F})$ be the family of $\mathcal{F}$-free graphs with ex$(n,\mathcal{F})$ edges. In particular, if $\mathcal{F}$ contains the unique  element, i.e., $\mathcal{F}=\{F\}$, then we denote by ex$(n,F)$ and EX$(n,F)$, for short.
Let $\Gamma_k$ be the family of graphs containing $k$-edge-disjoint triangles.
Gy\H{o}ri \cite{GE} determined the Tur\'{a}n number ex$(n,\Gamma_k)$ as follows.
\begin{thm}(\cite{GE})\label{thm-1}
Let $G$ be a graph of order $n$ that does not contain a subgraph belonging to $ \Gamma_k$, $k\geq 1$.
Then $e(G)\leq ex(n,\Gamma_k)=\lfloor\frac{n^2}{4}\rfloor+k-1$, and the extremal graph is obtained from a complete bipartite graph with color classes of size $\lceil \frac{n}{2}\rceil$ and $\lfloor \frac{n}{2}\rfloor$ by embedding edges of size $k-1$.
\end{thm}
As usual, given a graph $H$, let SPEX$(n,H)$ be the family of $H$-free graphs with the maximum spectral radius.
A $k$-fan, denoted by $F_k$, is the graph obtained from $k$ triangles by sharing a common vertex. Cioab\u{a}, Feng, Tait and Zhang \cite{CFTZ} proved that SPEX$(n,F_k)\subseteq \mbox{EX}(n,F_k)$.
To generalize their result, Li and Peng \cite{LP} showed that
SPEX$(n,H_{s,k})\subseteq \mbox{EX}(n,H_{s,k})$
where $H_{s,k}$ is the graph obtained from $s$ triangles and $k$ odd cycles of lengths at least 5 by sharing a common vertex.
The odd wheel $W_{2k+1}$ is the graph formed by joining a vertex to a cycle of length $2k$.
Very recently, Cioab\u{a}, Desai and Tait \cite{CDT} showed that $\mbox{SPEX}(n,W_5)\subseteq \mbox{EX}(n,W_5)$ and SPEX$(n,W_{2k+1})$ for $k\geq 3$ and $k\notin\{4,5\}$ is obtained from a complete bipartite graph with parts $L$ and $R$ of size $n+s$ and $n-s$ with $|s|\leq1$ by embedding a $(k-1)$-regular or nearly $(k-1)$-regular in $G[L]$ and exactly one edge in $G[R]$.
They posed the following interesting conjecture for further research.
\begin{conj}\cite{CDT}
Let $F$ be any graph such that the graphs in EX$(n, F)$ are Tur\'{a}n graphs plus $O(1)$ edges. Then SPEX$(n, F)\subseteq \mbox{EX}(n, F)$ for $n$ large enough.
\end{conj}

In this paper, we give a spectral version of Theorem \ref{thm-1},
in which the extremal graph is completely characterized.
This also provides a support for Cioab\u{a}, Desai and Tait's conjecture.

\begin{thm}\label{thm-2}
Let $k\geq 2$ and $G$ be a $\Gamma_k$-free graph of order $n$ sufficiently large.
If $G$ attains the maximum spectral radius, then
\begin{equation*}\label{th2-0}
  G\in \mbox{EX}(n, \Gamma_k).
\end{equation*}
More precisely, $G$ is obtained from $K_{\lfloor\frac{n}{2}\rfloor,\lceil\frac{n}{2}\rceil}$ by embedding a graph $H$ of size $k-1$
in the $\lfloor\frac{n}{2}\rfloor$-vertex partition set,
where $H\cong C_3$ for $k=4$ and $H\cong K_{1,k-1}$ otherwise.
\end{thm}

\section{Spectral condition on two cycles with the same length}\label{sec2}

Before beginning our proof, we first give some notation not defined in the above section. Let $G$ be a graph and  $d_G(v_i)$ be  the degree of the vertex $v_i$ in $G$. For $u,v \in V (G)$, denote  by $d_G(u, v)$ the distance between $u$ and $v$, i.e. the length of a shortest path between $u$ and $v$. Set $N_G^d(u)=\{v \mid v \in V (G), d_G(v, u) = d\}$. Specially, we use $N_G(u)$ instead of $N_G^1(u)$ and $N_G[u]=\{v \mid v \in N_G(u)\}\cup \{u\}$. Denote by $E_G(V_1,V_2)$ the set of edges between $V_1$ and $V_2$, and denoted $|E_G(V_1,V_2)|$ by $e_G(V_1,V_2)$. For the sake of simplicity, we shall omit all the subscripts if $G$ is clear from the context.

Denote by $G_{uv}$ the graph obtained from $G$ by subdividing the edge $uv$, that is, introducing a new vertex on the edge $uv$. Let $Y_n$ be the graph obtained from a path $v_1v_2\cdots v_{n-4}$ by attaching two pendant vertices to $v_1$ and two pendant vertices to $v_{n-4}$.
Hoffman and Smith \cite{HS75} proved the following result, which is an important tool in spectral graph theory.
\begin{lemma}{\rm(\cite{HS75})}\label{Le-HS}
Let $G$ be a connected graph with $uv \in E(G)$. If $uv$
belongs to an internal path of $G$ and
$G\ncong Y_n$, then $\rho(G_{uv}) < \rho(G)$.
\end{lemma}

Now we are in a position to give the proof of Theorem \ref{thm1.2}.

\noindent{\bf{Proof of Theorem \ref{thm1.2}}.}
Let $G^\star$ be a graph with the maximum spectral radius among all graphs without two cycles with the same length. Then $G^\star$ is connected, otherwise we add cut edges between components, the resulting graph will increase the spectral radius but will not produce a new cycle, which is impossible.
Now let $Z=(z_1, z_2, \ldots, z_n)^t$ be the Perron vector of $A(G^\star)$ and $z_{u^\star}=\max\{z_i: 1\leq i\leq n\}$.
Note that $K_{1, n-1}^+$ does not contain two cycles with the same length and
$X=\left(\frac{\sqrt2}{2}, \frac{1}{\sqrt{2(n-1)}},\ldots,\frac{1}{\sqrt{2(n-1)}}\right)^t$ is the Perron vector of $K_{1, n-1}$.
Obviously, $X\neq Z$.
It follows that $\rho(G^\star)\geq \rho(K_{1, n-1}^+)>X^tA(K_{1, n-1}^+)X=\sqrt{n-1}+\frac{1}{n-1},$
which is equivalent to
\begin{equation}\label{1}
  \rho^2(G^{\star})-\frac{2}{n-1}\rho(G^\star)+\frac{1}{(n-1)^2}-n+1>0.
\end{equation}
We first give the following claim which indicates that except $u^{\star}$ there is no cut vertex.
\begin{claim}\label{Claim1}
For each vertex $u\in V(G^{\star})\backslash \{u^{\star}\}$, $u$ is not a cut vertex.
\end{claim}
\renewcommand\proofname{\bf Proof}
\begin{proof}
By the way of contradiction, assume that $u$ is a cut vertex.
Then there exist two components of $G^\star\backslash\{u\}$, say $G_1$ and $G_2$. Without loss of generality, suppose that $u^{\star}\in V(G_1)$. Let $G'$ be the graph obtained by deleting the edges between $u$ and $N_{V(G_2)}(u)$ and adding the edges between $u^{\star}$ and $N_{V(G_2)}(u)$. The resulting graph still does not have two edge-disjoint cycles.
It follows that,
\begin{equation}\label{8}
  \rho(G')-\rho(G^{\star})\geq Z^t(A(G')-A(G^{\star}))Z=(z_{u^{\star}}-z_{u})\sum_{v\in N_{V(G_2)}(u)}z_{v}\geq 0.
\end{equation}
If $\rho(G')=\rho(G^{\star})$, then $z_{u^{\star}}=z_{u}$ and $Z$ is also the Perron vector of $A(G')$. But on the other hand, $$\rho(G')z_{u^{\star}}=\sum_{v\in N(u^{\star})}z_{v}+\sum_{v\in N_{V(G_2)}(u)}z_{v}>\sum_{v\in N(u^{\star})}z_{v}=\rho(G^{\star})z_{u^{\star}}.$$ It follows that $\rho(G')>\rho(G^{\star})$, which contradicts the fact that $G^{\star}$ has maximum spectral radius.
\end{proof}
Then, we have the following claim.
\begin{claim}\label{claim-4.1}
For each vertex $u\in V(G^\star)$, $e(N(u))\leq 1$ and $e(N(u), N^2(u))\leq |N^2(u)|+1$.
\end{claim}
\renewcommand\proofname{\bf Proof}
\begin{proof}
By the way of contradiction, we assume that either $e(N(u))\geq2$ or $e(N(u), N^2(u))\geq |N^2(u)|+2$. Then $G^\star$ contains two triangles or two $C_4$'s, a contradiction.
Then the result follows.
\end{proof}
Set $|N(u^\star)|=q$, $|N^2(u^\star)|=t$ $(t>0)$ and $N(u^\star)=\{v_{11}, v_{12}, \ldots, v_{1q}\}$, $N^2(u^\star)=\{v_{21}, v_{22}, \ldots, v_{2t}\}$.
For simplicity, we use $d_{1i}$ $(1\leq i\leq q)$ and $d_{2j}$ $(1\leq j\leq t)$ to denote $d_{N(u^\star)}(v_{1i})$ and $d_{N(u^\star)}(v_{2j})$, respectively. Now let $(1)\times z_{u^\star}$. Then
\begin{eqnarray}\label{2}
\begin{aligned}
&\rho^2(G^\star)z_{u^\star}-\frac{2}{n-1}\rho(G^\star)z_{u^\star}+(\frac{1}{(n-1)^2}-n+1)z_{u^\star} \\
&=d(u^\star)z_{u^\star}+\sum_{i=1}^qd_{1i}z_{1i}+\sum_{j=1}^td_{2j}z_{2j}-\frac{2}{n-1}\rho(G^\star)z_{u^\star}+(\frac{1}{(n-1)^2}-n+1)z_{u^\star},
\end{aligned}
\end{eqnarray}
it is easy to see that $(\ref{2})>0$ and
\begin{eqnarray*}
\begin{aligned}
(\ref{2})&\leq z_{u^\star}(d(u^\star)+\sum_{i=1}^qd_{1i}+\sum_{j=1}^td_{2j}-\frac{2}{n-1}\rho(G^\star)+\frac{1}{(n-1)^2}-n+1)\\
&= z_{u^\star}(d(u^\star)+2e(N(u^\star))+e(N(u^\star), N^2(u^\star))-\frac{2}{n-1}\rho(G^\star)+\frac{1}{(n-1)^2}-n+1).
\end{aligned}
\end{eqnarray*}
It follows that
\begin{eqnarray}\label{3}
d(u^\star)+2e(N(u^\star))+e(N(u^\star), N^2(u^\star))-\frac{2}{n-1}\rho(G^\star)+\frac{1}{(n-1)^2}-n+1>0.
\end{eqnarray}
Let $C(u^\star)=V(G^\star)\backslash\{N[u^\star]\cup N^2(u^\star)\}$. By Claim \ref{claim-4.1}, we have
$$e(N(u^\star), N^2(u^\star))\leq |N^2(u^\star)|+1=n-d(u^\star)-|C(u^\star)|.$$
Then $(\ref{3})$ becomes
\begin{eqnarray}\label{6}
2e(N(u^\star))-|C(u^\star)|-\frac{2}{n-1}\rho(G^\star)+\frac{1}{(n-1)^2}+1>0.
\end{eqnarray}
Since $\rho(G^\star)>\sqrt{n-1}+\frac{1}{n-1}$, by inequality $(\ref{6})$, we have $|C(u^\star)|=0$ if $e(N(u^\star))=0$ and $|C(u^\star)|\leq 2$ if $e(N(u^\star))=1$.

Note that $\rho(G^\star)z_{v_{2j}}=\sum\limits_{u\sim v_{2j}}z_u\leq d(v_{2j})z_{u^\star}$
and $\rho(G^\star)>\sqrt{n-1}\geq 5$ as $n\geq 26$. Then
\begin{align}\label{4}
(\ref{2})\leq&z_{u^\star}(d(u^\star)+\sum_{i=1}^qd_{1i}+\sum_{j=1}^td_{2j}\frac{d(v_{2j})}{\rho(G^\star)}-\frac{2}{n-1}\rho(G^\star)+\frac{1}{(n-1)^2}-n+1)\nonumber\\
=& z_{u^\star}(d(u^\star)+\sum_{i=1}^qd_{1i}+e(N(u^\star), N^2(u^\star))-\alpha-\frac{2}{n-1}\rho(G^\star)+\frac{1}{(n-1)^2}-n+1),
\end{align}
where $\alpha=e(N(u^\star), N^2(u^\star))-\sum\limits_{j=1}^td_{2j}\frac{d(v_{2j})}{\rho(G^\star)}$, that is, $\alpha=\sum\limits_{j=1}^td_{2j}-\sum\limits_{j=1}^td_{2j}\frac{d(v_{2j})}{\rho(G^\star)}$.
Before proceeding, we need the following fact.
\begin{fact}\label{fact1}
$5t-5\alpha<\sum_{j=1}^td(v_{2j})\leq t+1+2e(N^2(u^\star))+e(N^2(u^\star), N^3(u^\star)).$
\end{fact}
\renewcommand\proofname{\bf Proof}
\begin{proof}
Recall that there is at most one vertex $v_{2j}\in N^2(u^\star)$ such that $d_{2j}=2$, that is, $t\leq \sum\limits_{j=1}^td_{2j}\leq t+1$.
It implies that $\sum\limits_{j=1}^td(v_{2j})\leq t+1+2e(N^2(u^\star))+e(N^2(u^\star), N^3(u^\star)).$ In the following, we only need to show that $\sum\limits_{j=1}^td(v_{2j})>5t-5\alpha.$
If $\sum\limits_{j=1}^td_{2j}=t$, then $\alpha=t-\sum\limits_{j=1}^t\frac{d(v_{2j})}{\rho(G^\star)}$. If $\sum_{j=1}^td_{2j}=t+1$, without loss of generality, let $d_{2t}=2$, then $d_{N^2(u^\star)}(v_{2t})\leq 1$ (otherwise, one can find two $C_5$'s). Combining this with $|C(u^\star)|\leq 2$, we have $d(v_{2t})\leq 5$. Then $$\alpha=t-1-\sum\limits_{j=1}^{t-1}\frac{d(v_{2j})}{\rho(G^\star)}+2-2\frac{d(v_{2t})}{\rho(G^\star)}\geq
t-1-\sum\limits_{j=1}^{t-1}\frac{d(v_{2j})}{\rho(G^\star)}+1-\frac{d(v_{2t})}{\rho(G^\star)}=
t-\sum\limits_{j=1}^t\frac{d(v_{2j})}{\rho(G^\star)}.$$ Thus $\alpha\geq t-\sum\limits_{j=1}^t\frac{d(v_{2j})}{\rho(G^\star)}>t-\frac{1}{5}\sum\limits_{j=1}^td(v_{2j})$, as required.
\end{proof}
\begin{claim}\label{claim-4.2}
$e(N(u^\star))=1$.
\end{claim}
\renewcommand\proofname{\bf Proof}
\begin{proof}
By the way of contradiction, we assume that $e(N(u^\star))=0$.
Recall that the only possible cut vertex is $u^\star$. Moreover, from $(\ref{6})$, it follows that $|C(u^\star)|=0$.
Now we derive the proof by the following three cases.

\noindent{\bf{$\underline{\mbox{Case 1. }}$}} $e(N^2(u^\star))\geq 3$.
Without loss of generality, suppose that $e_1, e_2, e_3\in E(G^\star[N^2(u^\star)])$.
It is clear that each edge in $\{e_1, e_2, e_3\}$ is contained in either a $C_3$ or a $C_5$,
which implies that $G^\star$ contains two $C_3$'s or two $C_5$'s, a contradiction.

\noindent{\bf{$\underline{\mbox{Case 2.}}$}} $1\leq e(N^2(u^\star))\leq 2$.
Recall that $2e(N(u^\star))=\sum_{i=1}^qd_{1i}=0$ and $|C(u^\star)|=0$. Then $e(N^2(u^\star), N^3(u^\star))=0$ and
$$(\ref{4})\leq z_{u^\star}(-\frac{2}{n-1}\rho(G^\star)+\frac{1}{(n-1)^2}+1-\alpha).$$
From $0<(\ref{2})\leq (\ref{4})$, we have $-\frac{2}{n-1}\rho(G^\star)+\frac{1}{(n-1)^2}+1-\alpha>0$, which implies that $\alpha<1$.
Then by Fact \ref{fact1}, we have
$$5t-5<\sum_{j=1}^td(v_{2j})\leq t+1+2e(N^2(u^\star)),$$ which is equivalent to $2|N^2(u^\star)|<e(N^2(u^\star))+3\leq 5$, that is, $e(N^2(u^\star))=1$ and $|N^2(u^\star)|=2$. Then $2|N^2(u^\star)|<e(N^2(u^\star))+3=4$, which contradicts the fact that $|N^2(u^\star)|=2$.

\noindent{\bf{$\underline{\mbox{Case 3. }}$}} $e(N^2(u^\star))=0$.
Combining $e(N(u^\star), N^2(u^\star))\leq t+1$ with the fact that no vertex in $N(u^\star)$ is a cut vertex,
we have $t=|N^2(u^\star)|=1$, and then $G^\star\cong C_4\bullet K_{1,n-4}$. By Lemma \ref{Le-HS}, $\rho(C_4\bullet K_{1,n-4})<\rho(K_{1,n-2}^+)<\rho(K_{1,n-1}^+)$, a contradiction.
\end{proof}
 For each edge in $G^\star[V(G^\star)\backslash N[u^\star]]$, the two end-vertices of which have no common neighbor since $e(N(u^\star))=1$. It follows that $e$ is contained in a $5$-cycle for each edge $e\in G^\star[N^2(u^\star)]$, which implies that $e(N^2(u^\star))\leq 1$. The following claim gives the clear characterization of $G^\star$.
\begin{claim}\label{claim-4.3}
$G^\star\cong K_{1,n-1}^+$.
\end{claim}
\renewcommand\proofname{\bf Proof}
\begin{proof}
Recall that $|C(u^\star)|\leq 2$ and there exists no cut vertex in $V(G^\star)\backslash \{u^\star\}$. It follows that $|N^4(u^\star)|=0$ and $|N^3(u^\star)|\leq 2$.
Note that $\rho(G)\cdot z_u=\sum_{v\sim u}z_v\leq d(u)z_{u^\star}$ and $\rho(G)>\sqrt{n-1}\geq 5$ as $n\geq 26$. Then $z_u<\frac{d(u)}{5}z_{u^\star}$ for each $u\in V(G^\star)\backslash \{u^\star\}$.
We derive the proof by the following three cases.

\noindent{\bf{$\underline{\mbox{Case 1. }}$}} $|N^3(u^\star)|=2.$ Let $N^3(u^\star)=\{u_1,u_2\}$. Recall that $e(N(u^\star), N^2(u^\star))\leq n-d(u^\star)-|C(u^\star)|=t+1$. We assert $e(N(u^\star), N^2(u^\star))=t+1$. If not, then $e(N(u^\star), N^2(u^\star))=t=n-d(u^\star)-|C(u^\star)|-1$ and $(\ref{3})$ becomes $-\frac{2}{n-1}\rho(G^\star)+\frac{1}{(n-1)^2}<0$, a contradiction. It implies that $u^\star$ is contained in a $4$-cycle. Therefore, for each $u\in N^3(u^\star)$, if $w, v\in N^2(u^\star)$ are two neighbors of $u$, then $N_{N(u^\star)}(w)\neq N_{N(u^\star)}(v)$. It follows that $e(N^2(u^\star), N^3(u^\star))\leq 3$ (otherwise, there are two $C_6$'s) and $u_1u_2\in E(G^\star)$ since $N^2(u^\star)$ does not have cut vertex. By the above discussion, we have $\sum_{i=1}^qd_{1i}=2e(N(u^\star))=2$, $e(N^2(u^\star), N(u^\star))=t+1$, $e(N^2(u^\star))\leq 1$ and $e(N^2(u^\star), N^3(u^\star)\leq 3$. Then $$(\ref{4})=z_{u^\star}(-\frac{2}{n-1}\rho(G^\star)+\frac{1}{(n-1)^2}+1-\alpha).$$ Combining the fact with $0<(\ref{2})\leq (\ref{4})$, we have $\alpha<1$. By Fact \ref{fact1}, we have $5t-5<t+6$, that is, $t\leq 2$. Then $t=2$ since $N^2(u^\star)$ does not contain cut vertex and $e(N^2(u^\star))\neq 1$ since $u^\star$ belongs to a 4-cycles. Since $u_1$, $u_2$ have no common neighbor in $N^2(u^\star)$, so $e(N^2(u^\star), N^3(u^\star))=2$. Again by Fact \ref{fact1}, we have $5t-5<t+3$, that is, $t\leq 1$, which contradicts $t=2$.

\noindent{\bf{$\underline{\mbox{Case 2. }}$}} $|N^3(u^\star)|=1$. Let $N^3(u^\star)=\{u\}$. If $d(u)\geq 3$, then  $G^\star$ contains either two $C_4$'s or two $C_6$'s, which is impossible. Note that $u^\star$ is the unique possible cut vertex. Thus $d(u)=2.$ Without loss of generality, set $N(u)=\{v_{21},v_{22}\}$. Then $v_{21}v_{22}\notin E(G^\star)$ and $d_{21}=d_{22}=1$ (otherwise, either $G^\star$ contains two $C_4$'s or two $C_6$'s). So $d(v_{21})+d(v_{22})\leq 5$ since $e(N^2(u^\star))\leq 1$. Let $G_1=G^\star-\{uv_{21},uv_{22}\}+\{u^\star u\}$. Recall that $z_{v_{2i}}<\frac{d(u)}{5}z_{u^\star}$ for $i=1,2$. Then $$Z^T(\rho(G_1)-\rho(G^\star))Z=z_{u^\star}z_u-z_u(z_{v_{21}}+z_{v_{22}})> z_{u^\star}z_u(1-\frac{d(v_{21})+d(v_{22})}{5})\geq 0,$$ which implies that $\rho(G_1)>\rho(G^\star)$, a contradiction.

\noindent{\bf{$\underline{\mbox{Case 3. }}$}} $|N^3(u^\star)|=0$. Note that $e(N^2(u^\star))\leq 1$ and $e(N(u^\star),N^2(u^\star))\leq t+1$. Then $t\leq 3$ since no vertex in $N(u^\star)$ is a cut vertex. 
If $t=1$ or $t=3$, then there exists a vertex $\{v_{21}\}\in N^2(u^\star)$ such that $d(v_{21})=2$ and let $N(v_{21})=\{v_{11}, v_{12}\}\subseteq N(u^\star)$. Besides, $v_{11}v_{12}\notin E(G^\star)$. It is clear that $d(v_{11})+d(v_{12})\leq 5$ for $t=1$. Moreover, $d(v_{11})+d(v_{12})\leq 5$ also holds for $t=3$ since otherwise there are two $C_5$'s. Let $G_2=G^\star-\{v_{21}v_{11},v_{21}v_{12}\}+\{u^\star v_{21}\}$. By a similar proof to Case 2, we have $\rho(G_2)>\rho(G^\star)$, a contradiction. If $t=2$, say $\{v_{21},v_{22}\}\in N^2(u^\star)$, according to the fact that no vertex in $N(u^\star)$ is a cut vertex and $v_{21},v_{22}$ have no common neighbor, then $v_{21}v_{22}\in E(G^\star)$. It follows that $d_{21}=d_{22}=1$ (otherwise, there are two $C_5$'s), that is, $d(v_{21})=d(v_{22})=2$. Let $N_{N(u^\star)}(v_{21})=v_{11}$. Then $d(v_{11})\leq 3$ and $d(v_{11})+d(v_{22})\leq 5$. Let $G_3=G^\star-\{v_{21}v_{11},v_{21}v_{22}\}+\{u^\star v_{21}\}$. Similar to the proof of Case 2, we have $\rho(G_3)>\rho(G^\star)$, a contradiction. Therefore, $t=0$, i.e., $G^\star\cong K_{1, n-1}^+$.

Therefore, the proof of Theorem \ref{thm1.2} is complete.
\end{proof}

\section{Spectral conditions for edge-disjoint cycles}\label{sec1}

We first need the following well-known result as our lemma.
\begin{lemma}(Brouwer and Haemers \cite[p. 30]{BH}; Godsil and Royle \cite[pp. 196--198]{GR}.)\label{quotient}
Let $A$ be a real symmetric matrix, and let $B$ be an equitable quotient matrix of $A$. Then the eigenvalues of $B$ are also eigenvalues of $A$. Furthermore, if $A$ is nonnegative and irreducible, then $$\rho(A) = \rho(B).$$
\end{lemma}
%
Now, we shall give the proof of Theorem \ref{Thm:edgedisjointcycle}.
\subsection{Proof of Theorem \ref{Thm:edgedisjointcycle}}
Suppose that $G^{\star}$ is a graph attaining the maximum spectral radius among all graphs of order $n$ without two edge-disjoint cycles. Clearly, $G^{\star}$ is connected. Since the graph $K_4\bullet K_{1,n-4}$ does not contain two edge-disjoint cycles, we have $\rho(G^{\star})\geq \rho(K_4\bullet K_{1,n-4})$.

 Let $X=(\frac{\sqrt2}{2}, \frac{1}{\sqrt{2(n-1)}},\ldots,\frac{1}{\sqrt{2(n-1)}})^t$ be the Perron vector of $K_{1,n-1}$. Notice that $X$ cannot be the Perron vector of   $K_4\bullet K_{1,n-4}$. Then, by the  Rayleigh quotient,  we obtain
$$\rho(G^{\star})\geq \rho(K_4\bullet K_{1,n-4})>X'A(K_4\bullet K_{1,n-4})X=\sqrt{n-1}+\frac{3}{n-1},$$ which implies that
\begin{equation}\label{7}
  \rho^2(G^{\star})>n-1+\frac{6}{\sqrt{n-1}}+\frac{9}{(n-1)^2}.
\end{equation}

Let $Y=(y_{v_1}, y_{v_2}, \ldots, y_{v_n})^t$ be the unique positive eigenvector of  $\rho(G^\star)$ with  $\max\{y_{v_i}: 1\leq i\leq n\}=1$, and let $u^{\star}\in V(G^\star)$ be such that  $y_{u^{\star}}=1$. Similar to the proof of Theorem \ref{thm1.2}, $u^\star$ is the only possible cut vertex.  

\begin{claim}\label{Cl2}
  Aside from some triangles, there are no internal paths in $G^{\star}$.
\end{claim}
\renewcommand\proofname{\bf Proof}
\begin{proof}
By Lemma \ref{Le-HS}, contracting an edge in an internal path will increase the spectral radius but not increase the number of edge-disjoint cycles. Then the claim follows.
\end{proof}

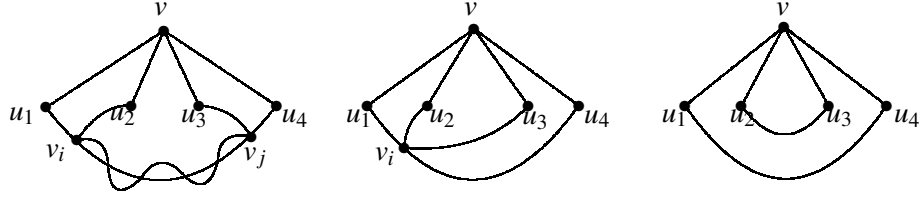
\begin{figure}[H]
 \begin{center}
\begin{picture}(331.3,94.9)(0,20)
\put(58.0,84.1){\circle*{4}}
\put(13.8,55.1){\circle*{4}}
\qbezier(58.0,84.1)(35.9,69.6)(13.8,55.1)
\put(45.7,55.8){\circle*{4}}
\qbezier(58.0,84.1)(51.8,70.0)(45.7,55.8)
\put(71.1,55.8){\circle*{4}}
\qbezier(58.0,84.1)(64.5,70.0)(71.1,55.8)
\put(100.1,55.8){\circle*{4}}
\qbezier(58.0,84.1)(79.0,70.0)(100.1,55.8)
\qbezier(13.8,55.1)(56.6,0.0)(100.1,55.8)
\put(25.4,42.8){\circle*{4}}
\qbezier(45.7,55.8)(34.8,56.6)(25.4,42.8)
\put(90.6,44.2){\circle*{4}}
\qbezier(71.1,55.8)(84.1,55.1)(90.6,44.2)
\qbezier(25.4,42.8)(36.3,45.0)(37.0,34.1)
\qbezier(37.0,34.1)(37.7,17.4)(49.3,28.3)
\qbezier(49.3,28.3)(56.6,39.9)(66.0,29.0)
\qbezier(66.0,29.0)(74.7,21.0)(78.3,35.5)
\qbezier(78.3,35.5)(78.3,46.4)(90.6,44.2)
\put(174.0,84.8){\circle*{4}}
\put(134.1,55.8){\circle*{4}}
\qbezier(174.0,84.8)(154.1,70.3)(134.1,55.8)
\put(156.6,55.8){\circle*{4}}
\qbezier(174.0,84.8)(165.3,70.3)(156.6,55.8)
\put(194.3,55.8){\circle*{4}}
\qbezier(174.0,84.8)(184.2,70.3)(194.3,55.8)
\put(213.2,55.8){\circle*{4}}
\qbezier(174.0,84.8)(193.6,70.3)(213.2,55.8)
\qbezier(134.1,55.8)(173.3,0.0)(213.2,55.8)
\put(147.9,39.9){\circle*{4}}
\qbezier(156.6,55.8)(148.6,50.8)(147.9,39.9)
\qbezier(147.9,39.9)(177.6,37.7)(194.3,55.8)
\put(290.0,85.6){\circle*{4}}
\put(253.0,55.8){\circle*{4}}
\qbezier(290.0,85.6)(271.5,70.7)(253.0,55.8)
\put(274.1,55.8){\circle*{4}}
\qbezier(290.0,85.6)(282.0,70.7)(274.1,55.8)
\put(306.7,55.8){\circle*{4}}
\qbezier(290.0,85.6)(298.3,70.7)(306.7,55.8)
\put(328.4,55.8){\circle*{4}}
\qbezier(290.0,85.6)(309.2,70.7)(328.4,55.8)
\qbezier(253.0,55.8)(290.7,0.7)(328.4,55.8)
\qbezier(274.1,55.8)(290.0,34.1)(306.7,55.8)
\put(170.4,94.9){\makebox(0,0)[tl]{$v$}}
\put(54.4,94.9){\makebox(0,0)[tl]{$v$}}
\put(287.1,94.9){\makebox(0,0)[tl]{$v$}}
\put(0.0,55.3){\makebox(0,0)[tl]{$u_1$}}
\put(13.8,41.3){\makebox(0,0)[tl]{$v_i$}}
\put(126.2,54.4){\makebox(0,0)[tl]{$u_1$}}
\put(244.3,54.6){\makebox(0,0)[tl]{$u_1$}}
\put(37.7,55.1){\makebox(0,0)[tl]{$u_2$}}
\put(157.3,54.4){\makebox(0,0)[tl]{$u_2$}}
\put(269.7,54.6){\makebox(0,0)[tl]{$u_2$}}
\put(137.0,39.9){\makebox(0,0)[tl]{$v_i$}}
\put(192.1,53.7){\makebox(0,0)[tl]{$u_3$}}
\put(306.0,54.6){\makebox(0,0)[tl]{$u_3$}}
\put(64.5,54.4){\makebox(0,0)[tl]{$u_3$}}
\put(87.7,40.2){\makebox(0,0)[tl]{$v_j$}}
\put(331.3,54.3){\makebox(0,0)[tl]{$u_4$}}
\put(102.2,54.4){\makebox(0,0)[tl]{$u_4$}}
\put(215.3,54.1){\makebox(0,0)[tl]{$u_4$}}
\end{picture}
\caption{All cases of $(u_1,u_4)$-path and $(u_2,u_3)$-path.}
\label{fig1}
\end{center}
\end{figure}

\begin{claim}\label{Cla4}
  Each  vertex of $V(G^{\star})\backslash \{u^{\star}\}$ has degree at most  $3$.
\end{claim}
\renewcommand\proofname{\bf Proof}
\begin{proof}
By contradiction, assume that there exists some $v\in V(G^{\star})\backslash \{u^{\star}\}$ such that $d(v)\geq 4$.
Let $u_1,u_2,u_3,u_4$ be four neighbors of $v$.
By Claim \ref{Claim1}, $v$ is not a cut vertex.
Thus, there is a path $P_1$ from $u_1$ to $u_4$ with $v\notin V(P_1)$.
We call it a $(u_1,u_4)$-path.
Similarly, we can find a $(u_2,u_3)$-path $P_2$ with $v\notin V(P_2)$.
Let $|V(P_1)\cap V(P_2)|=k$.
Whether $k\geq2$, $k=1$ or $k=0$ (\emph{see Figure. \ref{fig1}}),
we can always find two edge-disjoint cycles, a contradiction.
\end{proof}

Suppose that $N(u^{\star})=\{v_{11}, v_{12}, \ldots, v_{1q}\}$ and  $N^2(u^{\star})=\{v_{21}, v_{22}, \ldots, v_{2t}\}$. For simplicity, let  $d_{1i}=d_{G^{\star}[N(u^{\star})]}(v_{1i})$  $(1\leq i\leq q)$  and $d_{2j}=d_{G^{\star}[N(u^{\star})]}(v_{2j})$ $(1\leq j\leq t)$.
%
 Our goal of the rest is to show that $e(N(u^{\star}))=3$.
We first show the following claim to insure that $N(u^{\star})$ is not an independent set.
\begin{claim}\label{Cla5}
   $e(N(u^{\star}))\geq1$
\end{claim}
\renewcommand\proofname{\bf Proof}
\begin{proof}
Suppose to the contrary that $e(G^\star[N(u^{\star})])=0$. Let $B_1, B_2,\ldots, B_l$ be the components of  $V(G^{\star})\backslash N[u^{\star}]$. Clearly,   $l\geq 1$ because $G^{\star}\not\cong K_{1,n-1}$ due to $\rho(G^{\star})>\rho(K_{1,n-1})$.
By Claim \ref{Claim1}, each vertex of $V(G^\star)\backslash\{u^{\star}\}$ is not a cut vertex.
Then $|N_{N(u^{\star})}(B_i)|\geq2$ and then $G^{\star}[N[u^{\star}]\cup B_i]$ contains a cycle for $1\leq i\leq l$. We derive the proof by the following two cases.

\noindent{\bf{$\underline{\mbox{Case 1.}}$}} $l\geq 2$.
We claim that $G^{\star}[B_i]$ is a tree for $1\leq i\leq l$.
If not, assume that $G^{\star}[B_j]$ contains a cycle, say $C_1$ for some $j\in [1,l]$.
Then $C_1$ together with a cycle in $G^{\star}[N[u^{\star}]\cup B_i]$ $(i\neq j)$ are two edge-disjoint cycles in $G^\star$, a contradiction.
Similarly, we have $e(v, B_i)\leq 1$ for each vertex $v\in N(u^{\star})$.
Then we claim that $|B_i|=1$ for $1\leq i\leq l$.
Otherwise, suppose that $|B_j|\geq 2$ for some $j\in [1,l]$, then there are at least two leaves in $B_j$,
say $v_{21}$ and $v_{22}$. Then $N_{N(u^{\star})}(v_{21})\cap N_{N(u^{\star})}(v_{22})=\varnothing$.
Moreover, by Claims \ref{Claim1} and \ref{Cl2}, we have $|N_{N(u^{\star})}(v_{21})|=|N_{N(u^{\star})}(v_{22})|=2$.
Without loss of generality, assume that $N_{N(u^{\star})}(v_{21})=\{v_{11}, v_{12}\}$ and $N_{N(u^{\star})}(v_{22})=\{v_{13}, v_{14}\}$.
Then $u^{\star}v_{11}v_{21}v_{12}u^{\star}$ and $u^{\star}v_{13}v_{22}v_{14}u^{\star}$ are two edge-disjoint cycles, a contradiction.
Then by Claim 2, $e(N(u^{\star}), B_i)=3$ for $1\leq i\leq l$. It follows that $l=2$ (otherwise, one may easily find two edge-disjoint cycles).
Then $G^{\star}\cong K_{3, 3}\bullet K_{1, n-6}$.
Notice that   $$B=\begin{bmatrix}0 & 1 & 0 & 0\\ n-6 & 0 & 0 & 3\\ 0 & 0 & 0 & 3\\ 0 & 1 & 2 & 0 \end{bmatrix}.$$
is an  equitable quotient matrix of $A(K_{3, 3}\bullet K_{1, n-6})$, and the  characteristic polynomial of $B$ is $f(x)=x^4-(n+3)x^2+6n-36.$
By Lemma \ref{quotient}, it follows that $$\rho^2(G^{\star})=\rho^2(B)=\frac{n+3+\sqrt{(n+3)^2-4(6n-36)}}{2}<n-1$$ as $n\geq 17$, a contradiction.

\noindent{\bf{$\underline{\mbox{Case 2.}}$}} $l=1$.
We claim that $|B_1|\geq 2$. Otherwise, suppose that $|B_1|=1$ and $B_1=\{v\}$, then the vertex $w\in N(v)$ lies in an internal path, a contradiction.
If $|B_1|=2$ then $G^{\star}\cong G_1$ (as shown in Figure. \ref{fig2}).
Let $G_2=G_1-\{v_{11}v_{21},v_{11}v_{22}\}+\{v_{21}u^{\star},v_{22}u^{\star}\}$ (see Figure. \ref{fig2}).
Note that $G_2$ does not contain two edge-disjoint cycles. However,
$$Y^t(\rho(G_2)-\rho(G^{\star}))Y\geq 2(y_{v_{21}}+y_{v_{22}})(y_{u^{\star}}-y_{v_{11}})>0$$ since $y_{v_{11}}\leq \frac{3y_{u^{\star}}}{\rho({G^{\star}})}<y_{u^{\star}}$ as $n\geq 17$.
It follows that $\rho(G^{\star})<\rho(G_2)$, a contradiction.
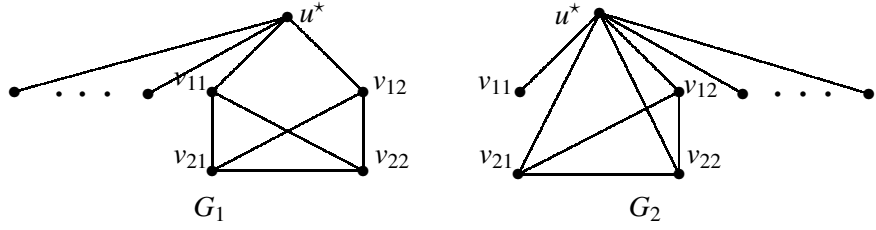
\begin{figure}[H]
\begin{center}
\begin{picture}(319.7,90.6)(0,10)
\put(102.2,84.8){\circle*{4}}
\put(74.0,56.6){\circle*{4}}
\qbezier(102.2,84.8)(88.1,70.7)(74.0,56.6)
\put(130.5,56.6){\circle*{4}}
\qbezier(102.2,84.8)(116.4,70.7)(130.5,56.6)
\put(74.0,26.8){\circle*{4}}
\qbezier(74.0,56.6)(74.0,41.7)(74.0,26.8)
\put(130.5,26.8){\circle*{4}}
\qbezier(130.5,56.6)(130.5,41.7)(130.5,26.8)
\qbezier(74.0,26.8)(102.2,26.8)(130.5,26.8)
\qbezier(74.0,56.6)(102.2,41.7)(130.5,26.8)
\qbezier(130.5,56.6)(102.2,41.7)(74.0,26.8)
\put(219.0,86.3){\circle*{4}}
\put(189.2,56.6){\circle*{4}}
\qbezier(219.0,86.3)(204.1,71.4)(189.2,56.6)
\put(188.5,25.4){\circle*{4}}
\qbezier(219.0,86.3)(203.7,55.8)(188.5,25.4)
\put(248.7,56.6){\circle*{4}}
\qbezier(219.0,86.3)(233.8,71.4)(248.7,56.6)
\put(248.7,25.4){\circle*{4}}
\qbezier(248.7,56.6)(248.7,41.0)(248.7,25.4)
\qbezier(188.5,25.4)(218.6,41.0)(248.7,56.6)
\qbezier(188.5,25.4)(218.6,25.4)(248.7,25.4)
\qbezier(219.0,86.3)(233.8,55.8)(248.7,25.4)
\put(106.6,90.6){\makebox(0,0)[tl]{$u^{\star}$}}
\put(202.3,89.9){\makebox(0,0)[tl]{$u^{\star}$}}
\put(59.5,63.8){\makebox(0,0)[tl]{$v_{11}$}}
\put(134.1,62.4){\makebox(0,0)[tl]{$v_{12}$}}
\put(173.5,62.4){\makebox(0,0)[tl]{$v_{11}$}}
\put(250.3,60.9){\makebox(0,0)[tl]{$v_{12}$}}
\put(59.5,34.3){\makebox(0,0)[tl]{$v_{21}$}}
\put(174.0,34.3){\makebox(0,0)[tl]{$v_{21}$}}
\put(134.9,34.3){\makebox(0,0)[tl]{$v_{22}$}}
\put(251.6,34.3){\makebox(0,0)[tl]{$v_{22}$}}
\put(34.8,55.8){\circle*{2}}
\put(25.4,55.8){\circle*{2}}
\put(16.7,55.8){\circle*{2}}
\put(0.0,56.6){\circle*{4}}
\qbezier(102.2,84.8)(51.1,70.7)(0.0,56.6)
\put(285.7,55.8){\circle*{2}}
\put(295.1,55.8){\circle*{2}}
\put(304.5,55.8){\circle*{2}}
\put(319.7,55.8){\circle*{4}}
\qbezier(219.0,86.3)(269.3,71.1)(319.7,55.8)
\put(67,17){\makebox(0,0)[tl]{$G_1$}}
\put(230,17){\makebox(0,0)[tl]{$G_2$}}
\put(272.6,55.8){\circle*{4}}
\qbezier(219.0,86.3)(245.8,71.1)(272.6,55.8)
\put(50.0,55.8){\circle*{4}}
\qbezier(102.2,84.8)(76.1,70.3)(50.0,55.8)
\end{picture}
\caption{The graphs $G_1$ and $G_2$.}
\label{fig2}
\end{center}
\end{figure}
So in the following, we may assume that $|B_1|\geq 3$. Note that $$\rho^2({G^{\star}})=\rho^2({G^{\star}})y_{u^{\star}}=d(u^{\star})+\sum_{i=1}^qd_{1i}y_{1i}+\sum_{j=1}^td_{2j}y_{2j}.$$
According to inequality $(\ref{7})$ and the assumption that $e(N(u^\star))=0$, we have
$$n-1+\frac{6}{\sqrt{n-1}}+\frac{9}{(n-1)^2}<d(u^{\star})+\sum_{j=1}^td_{2j}y_{2j}\leq d(u^{\star})+e(N(u^{\star}), B_1).$$
It follows that $e(N(u^{\star}), B_1)>n-1-d(u^{\star})=|B_1|$.
Then there exists some  $v_{21}\in B_1$ such that $d_{21}\geq 2$.
Thus $G^\star[B_1]$ is a tree.
Otherwise, a cycle in $G^\star[B_1]$ together with a cycle in $G^{\star}[v_{21}\cup N[u^{\star}]]$ are two edge-disjoint cycles, a contradiction.
By Claims \ref{Claim1} and \ref{Cla4}, we have $2 \leq d(v)\leq 3$ for each $v\in V(B_1)$.
Set $S=\{v\in V(B_1)| d(v)=2\}$. By Claim \ref{Cl2},
each vertex  of $S$ belongs to some triangle of $G^\star.$
We aim to show that $|S|=0$. By contradiction, suppose that $v_{22}\in S$, then $v_{21}v_{22}\notin E(G^{\star})$ since $|B_1|\geq 3$ and $d(v_{21})=3$.
Then the $4$-cycle from $G^\star[N[u^\star]\cup \{v_{21}\}]$ together with the triangle containing $v_{22}$ are two edge-disjoint cycles, a contradiction.
Therefore, $e(N(u^{\star}), B_1)=3|B_1|-2(|B_1|-1)=|B_1|+2\geq 5$.
Thus $|N_{N(u^{\star})}(B_1)|\geq 3$ by Claim \ref{Cla4}.
Suppose that $\{v_{11}, v_{12}, v_{13}\}\subseteq N(B_1)$ and $N_{N(u^{\star})}(v_{21})=\{v_{11}, v_{12}\}$.
Since $e(v_{13}, B_1)=2$ by Claims \ref{Cl2} and \ref{Cla4}, the cycle $v_{21}v_{11}u^{\star}v_{12}v_{21}$ together with a cycle in $G^\star[B_1\cup \{v_{13}\}]$ are two edge-disjoint cycles, a contradiction.
Therefore, the proof of Claim \ref{Cla5} is complete.
\end{proof}
According to Claim \ref{Cla5}, we see that $e(G^{\star}[N(u^{\star})])\geq1$. Now we shall prove that $V(G^\star)=N[u^\star]$.
If not, suppose that $B_1, B_2,\ldots, B_l$ are components of $G^\star-N[u^\star]$.
Then each $G^\star[B_i]$ is a tree and $e(v, B_i)\leq 1$ for each $v\in N(u^\star)$.
Furthermore,  we claim that $|B_i|=1$ for $i=1,\ldots,l$, since otherwise  there are two leaves in $G^\star[B_i]$, and we can find two edge-disjoint $C_4$'s in $G^\star$, a contradiction.
If $l\geq 2$, 
we set  $B_1=\{u\}$ and $B_2=\{v\}$.
If $|N(u)\cap N(v)|\geq 2$, then a triangle in $G^\star[N[u^\star]]$ together with a $C_4$ in $G^\star[N(u^\star)\cup \{u,v\}]$ are two edge-disjoint cycles, a contradiction.
If $|N(u)\cap N(v)|\leq 1$, then again we will find two edge-disjoint cycles since $d(u), d(v)\geq 3$, a contradiction.
It follows that $l=1$. Set $B_1=\{u\}$. Then $2\leq d(u)\leq 3$. If $d(u)=3$ and $N(u)=\{v_{11},v_{12},v_{13}\}$,
again by Claims \ref{Cl2} and \ref{Cla4}, we have $d(v_{1i})=3$ for $i=1,2,3$, then $G^\star$ contains two edge-disjoint triangles, a contradiction. If $d(u)=2$ and $N(u)=\{v_{11},v_{12}\}$, by Claim \ref{Cl2}, we have $v_{11}v_{12}\in E(G^\star)$, then $G^\star+u^\star u$ does not contain two edge-disjoint cycles but $\rho(G^\star+u^\star u)>\rho(G^\star)$, a contradiction.
Note that $e(N(u^\star))\leq 3$.
Then by the maximality of the spectral radius, we have $G^\star\cong K_4\bullet K_{1,n-4}$.

\subsection{Proof of Theorem \ref{thm-2}}
First we list some lemmas, which are useful in the proof of Theorem \ref{thm-2}.
\begin{lemma}(\cite{CFTZ})\label{lower bound with triangles}
Let $G$ be a graph, and let $t$ denote the number of triangles in $G$. Then $$e(G) \geq \rho^2(G) - \frac{3t}{\rho(G)}.$$
\end{lemma}

\begin{lemma}(Cioab\u{a}, Feng, Tait and Zhang, \cite{CFTZ})\label{thm-0}
Let $G$ be an $F_k$-free graph of order $n$. For sufficiently large $n$, if $G$ has the maximal spectral radius, then
$$G \in \textrm{EX}(n, F_k).$$
\end{lemma}

Denote by $\beta(G)$ and $\Delta(G)$ be the matching number and the maximum degree of $G$, respectively.
 For  any two positive integers $\beta$ and $\Delta$, we define $f(\beta, \Delta)=\max\{|E(G)|: \beta(G)\leq \beta, \Delta(G)\leq \Delta \}$. Chv\'atal and Hanson \cite{Chvatal76} proved the following result.

\begin{lemma}[Chv\'atal and Hanson \cite{Chvatal76}]\label{Chvatal76}
For positive integers $\beta$, $\Delta \geq 1$,
$$f(\beta, \Delta)= \Delta \beta +\left\lfloor\frac{\Delta}{2}\right\rfloor
 \left \lfloor \frac{\beta}{\lceil{\Delta}/{2}\rceil }\right \rfloor
 \leq \Delta \beta+\beta.$$
\end{lemma}

\begin{lemma}(\cite{WXH})\label{wu}
Let $G$ be a connected graph and let $X$ be the eigenvector corresponding to $\rho(G)$.
If $x_u\geq x_v$, and let $G'=G-\{vw|w\in N(v)\backslash N(u)\}+\{uw|w\in N(v)\backslash N(u)\},$ then $\rho(G) < \rho(G')$.
\end{lemma}

Denote by $K_{\lceil\frac{n}{2}\rceil,\lfloor\frac{n}{2}\rfloor}\diamond K_{1,k-1}$ the graph obtained  by embedding a copy of $K_{1,k-1}$ in the part of  order $\lfloor\frac{n}{2}\rfloor$ in  $K_{\lceil\frac{n}{2}\rceil,\lfloor\frac{n}{2}\rfloor}$.
Next, we give the  proof of Theorem \ref{thm-2}.
\vspace{2mm}

\noindent{\bf{Proof of Theorem \ref{thm-2}}}.
Assume that $G^\star$ is a graph attaining the maximum spectral radius among all graphs of order $n$ containing  no \emph{$\Gamma_k$}. Clearly, $G^\star$ is connected. Let us first prove that  $G^\star\in \textrm{EX}(n,\Gamma_k)$. Since  the proof method is almost the same as  the one in Lemma \ref{thm-0} (cf. \cite{CFTZ}), we omit the details.

\begin{claim}\label{claim-k}
$G^\star\in \textrm{EX}(n,\Gamma_k).$
\end{claim}
\renewcommand\proofname{\bf Proof}
\begin{proof}
Notice that $K_{\lceil\frac{n}{2}\rceil,\lfloor\frac{n}{2}\rfloor}\diamond K_{1,k-1}$ is $\Gamma_k$-free. By the maximality of $\rho(G^\star)$, we have
$\rho(G^\star)\geq \rho(K_{\lceil\frac{n}{2}\rceil,\lfloor\frac{n}{2}\rfloor}\diamond K_{1,k-1})\geq \frac{2}{n}e(K_{\lceil\frac{n}{2}\rceil,\lfloor\frac{n}{2}\rfloor}\diamond K_{1,k-1})=\frac{2}{n}(\lfloor\frac{n^2}{4}\rfloor+k-1)>\frac{n}{2}$. Also, since $G^\star$ is $\Gamma_k$-free, $F_k$ cannot be a subgraph of $G^\star$, i.e., for any $v\in V(G^\star)$, $G^\star[N(v)]$  contains no  $kK_2$.  Let $t$ be the number of triangles in $G^\star.$
By Lemma \ref{Chvatal76},
$$3t=\sum_{v\in V(G^\star)} e(G^\star[N(v)]) \leq \sum_{v\in V(G^\star)} \mathrm{ex}(d(v), kK_2) \leq \sum_{v\in V(G^\star)} \mathrm{ex}(n, kK_2) \leq \sum_{v\in V(G^\star)} kn=kn^2, $$
which gives that  $t\leq \frac{kn^2}{3}.$
Therefore, by Lemma \ref{lower bound with triangles},  we have $e(G^\star)\geq \rho^2(G^\star)-\frac{6t}{n}\geq \rho^2(G^\star)-2kn> \frac{n^2}{4}-2kn.$ Let $\varepsilon$ and $\delta$ be fixed positive constants with $\delta<\frac{1}{10(k+1)^2}$, $ \varepsilon<\frac{\delta^2}{16}$.
 Notice that $\textrm{ex}(n, \Gamma_k)= \lfloor\frac{n^2}{4}\rfloor+k-1\leq \textrm{ex}(n, F_k)\leq\lfloor\frac{n^2}{4}\rfloor+k^2-k$. Then, according to the proof in \cite[Lemma 10]{CFTZ}, we have
 \begin{fact}\label{maxcut}
 The graph  $G^\star$ has a partition $V(G^\star)=S\cup T$ which gives  a  maximum cut such that $e(S,T)\geq (\frac{1}{4}-\varepsilon)n^2.$
Furthermore, $\left(\frac{1}{2} - \sqrt{\varepsilon} \right) n \leq |S|, |T| \leq \left(\frac{1}{2} + \sqrt{\varepsilon}\right) n.$
\end{fact}
Let $L= \left\{v\in V(G^\star): d(v) \leq \left(\frac{1}{2}-\frac{1}{4(k+1)}\right) n\right\}$ and $W= \left\{ v\in S: d_S(v) \geq \delta n\right\} \cup \left\{v \in T: d_T(v) \geq \delta n\right\}.$
We assert that $|L|<16k^2$. If not, there exists some $L'\subseteq L$ with $|L'|=16k^2$. Then
 \begin{eqnarray*}
e(G^\star-L') \ge e(G^\star)-\sum_{v\in L^{\prime}}d(v)\ge \frac{n^2}{4}-2kn-16k^2 \left(\frac{1}{2}-\frac{1}{4(k+1)} \right)n> \frac{(n-16k^2)^2}{4}+k-1
     \end{eqnarray*}
for sufficiently large $n$. By Theorem \ref{thm-1}, $G^\star-L'$ contains $k$ edge-disjoint triangles, and so is $G^\star$, a contradiction. Combing this with Fact \ref{maxcut}, as in \cite[Lemmas 13--17]{CFTZ}, we can prove the following results successively:
\begin{itemize}
\item $|W| <\frac{2\varepsilon}{\delta} n+\frac{2k^2}{\delta n}$ and $W\subseteq L$;
\item $L=\emptyset$, and both $G^\star[S]$ and $G^\star[T]$ are $K_{1,k}$- and $kK_2$-free;
\item $e(G^\star)\ge \frac{n^2}{4}-12k^2$, $\frac{n}{2}-4k\leq |S|, |T|\leq \frac{n}{2}+4k,$
and
$\frac{n}{2}-14k^2\le \delta(G^\star)\le \rho(G^\star)\le \Delta (G^\star)\le \frac{n}{2}+5k$;
\item For any $u\in V(G^\star)$,   $\mathbf{x}_u\geq 1-\frac{120 k^2}{n}$;
\item  $\big||S| - |T|\big| \leq 1$.
\end{itemize}

Now we prove  that $G^\star\in \textrm{EX}(n,\Gamma_k)$. By contradiction, assume that  $e(G^\star)\leq \textrm{ex}(n,\Gamma_k)-1.$  By Theorem \ref{thm-1}, every graph in $\textrm{EX}(n, \Gamma_k)$ has a maximum cut $(S,T)$ of size $\lfloor n^2/4\rfloor$. Since $\big||S| - |T|\big| \leq 1$, there exists some $H\in \textrm{EX}(n,\Gamma_k)$ with $V(H)=V(G^\star)$ such that the edges between  $S$ and $T$  in $H$ form a complete bipartite graph.
Let $E_+=E(H)\setminus E(G^\star)$ and $E_-=E(G^\star)\setminus E(H)$. Then $(E(G^\star) \cup E_+) \setminus E_- = E(H)$, and  $|E_+| \geq |E_-| + 1$ because $|E(G^\star)\cap E(H)|+|E_-|=e(G^\star)< e(H)=|E(G^\star)\cap E(H)|+|E_+|$. Also note that  $|E_-| \leq e(S) + e(T) < 2k^2$. Thus,
\begin{align*}
\rho(H) &\geq \frac{\mathbf{x}^T A(H) \mathbf{x}}{\mathbf{x}^T\mathbf{x}} \\
&= \rho(G^\star) +\frac{2}{\mathbf{x}^T\mathbf{x}} \sum_{ij\in E_+} \mathbf{x}_i\mathbf{x}_j - \frac{2}{\mathbf{x}^T\mathbf{x}} \sum_{ij\in E_-} \mathbf{x}_i\mathbf{x}_j \\
& \geq \rho(G^\star) + \frac{2}{\mathbf{x}^T\mathbf{x}} \left( |E_+|\left(1 - \frac{120k^2}{n}\right)^2 - |E_-|\right)\\
& \geq \rho(G^\star) + \frac{2}{\mathbf{x}^T\mathbf{x}} \left( |E_+|  - |E_-|-\frac{240k^2}{n}  |E_+| +\frac{(120k^2)^2}{n^2}  |E_+| \right)\\
& \geq \rho(G^\star) + \frac{2}{\mathbf{x}^T\mathbf{x}} \left( 1-\frac{240k^2}{n}  |E_+| +\frac{(120k^2)^2}{n^2}  |E_+| \right)\\
&>\rho(G^\star)
\end{align*}
for sufficiently  large $n$, contrary to the maximality of  $\rho(G^\star)$. Therefore, we conclude that  $e(G^\star)=\textrm{ex}(n,\Gamma_k)$, and the result follows.
\end{proof}

By Claim \ref{claim-k} and  Theorem \ref{thm-1}, we see that $G^\star$ has a partition
$V(G^\star)=S\cup T$ with $|S|=s=\lfloor\frac{n}{2}\rfloor$ and $|T|=t=\lceil\frac{n}{2}\rceil$  such that the edges between $S$ and $T$ form a complete bipartite graph and  $e(S)+e(T)=k-1$. Now we focus on characterizing the exact structure  of $G^\star$.

Let $X=(x_v:v\in V(G^\star))^t$ be the Perron vector of $G^\star$. Set $S=\{u_0,u_1,\ldots,u_{s-1}\}$ and $T=\{v_0,v_1,\ldots,v_{t-1}\}$. Without loss of generality,  suppose that $x_{u_0}\geq x_{u_1}\geq \cdots\geq x_{u_{s-1}}$ and  $x_{v_0}\geq x_{v_1}\geq \cdots\geq x_{v_{t-1}}$. We first assert that
\begin{equation}\label{equ::add}
N_{S}(u_{j})\subseteq N_{S}[u_{i}] ~~\textrm{for}~~ 0\leq i < j \leq s-1,
~~ N_{T}(v_{j})\subseteq N_{T}[v_{i}] ~~\textrm{for}~~ 0\leq i< j \leq t-1.
\end{equation}
Without loss of generality, we only prove the former by the way of contradiction.
Suppose that there exist $i,j$ with $i<j$ such that $N_{S}(u_{j})\not\subseteq N_{S}[u_{i}]$. Let $u\in N_{S}(u_{j})\setminus N_{S}[u_{i}]$ and  $G'=G^\star-uu_j+uu_i$. By Lemma \ref{wu},  $\rho(G')>\rho(G^\star)$. This is a contradiction, and so  (\ref{equ::add}) follows.  According to (\ref{equ::add}), we have  $N_{S}(u_{j})\subseteq N_{S}[u_{0}]$ for all $j\in[1,s-1]$ (resp. $N_{T}(v_{j})\subseteq N_{T}[v_{0}]$ for all $j\in[1,t-1]$), and it follows that $E(S)=E(N_S[u_0])$ (resp. $E(T)=E(N_T[v_0])$). Similarly, we can deduce that $E(N_S(u_0))=E(N_S[u_1])$ and $E(N_T(v_0))=E(N_T[v_1])$.

In the remaining part of the proof, we will frequently construct a $\Gamma_k$-free graph $G'$  from $G^\star$ by modifying some edges.  For convenience, we always use  $Y=(y_v:v\in V(G'))^t$ to denote the Perron vector of $G'$. Also, we write $\rho'=\rho(G')$ and $\rho=\rho(G^\star)$ for short.

Let $s_0=d_S(u_0)$ and $t_0=d_T(v_0)$. Clearly, $e(S)\geq s_0$ and $e(T)\geq t_0$. In the following, we shall analyze the structure of $G^\star$ according to the structure of $G^\star[S]$ and $G^\star[T]$. Before proceeding, we need the following two claims.



%

\begin{claim}\label{claim-k2}
$e(S)\in \{s_0, s_0+1\}$ and $e(T)\in\{t_0,t_0+1\}$.
\end{claim}
\renewcommand\proofname{\bf Proof}
\begin{proof}
We first prove that $e(S)\in \{s_0, s_0+1\}$. By contradiction, suppose that  $e(S)\geq s_0+2$. Let $E^\star=E(N_S(u_0))$. Since $e(S)=e(N_S[u_0])\geq d_{S}(u_{0})+2$,  we have $|E^\star|\geq2$.
Let $E'=\{u_0u_i|s_0+1\leq i\leq e(S)\}$ and $G'=G^\star-E^\star+E'$.
Then $|E'|=e(S)-s_0=|E^\star|$, $G'[S]\cong K_{1,e(S)}\cup (s-e(S)-1)K_1$,   and $G'$ is $\Gamma_k$-free. Recall that $X$ and $Y$ are the Perron vectors of $G^\star$ and $G'$, respectively. Then
$$
Y^t(\rho'-\rho)X=Y^t(A(G')-A(G^\star))X=\sum_{i=s_0+1}^{e(S)}(x_{u_0}y_{u_i}+y_{u_0}x_{u_i})-\sum_{u_iu_j\in E^\star}(x_{u_i}y_{u_j}+x_{u_j}y_{u_i}).
$$
Observe that $x_{u_i}=x_{u_j}$ for $i,j\in[s_0+1,e(S)]$ and $y_{u_i}=y_{u_j}$ for $i,j\in[1, e(S)]$.  Then
$$
\sum_{i=s_0+1}^{e(S)}(x_{u_0}y_{u_i}+y_{u_0}x_{u_i})=|E^\star|(x_{u_0}y_{u_1}+y_{u_0}x_{u_{s_0+1}})
$$
and
$$
\sum_{u_iu_j\in E^\star}(x_{u_i}y_{u_j}+x_{u_j}y_{u_i})\leq |E^\star|y_{u_1}(x_{u_1}+x_{u_2})\leq 2|E^\star|x_{u_0}y_{u_1}.
$$
It follows that
$Y^t(\rho'-\rho)X\geq|E^\star|(y_{u_0}x_{u_{s_0+1}}-x_{u_0}y_{u_1}).$
Note that $\rho x_{u_0}-\rho x_{u_{s_0+1}}=\sum_{i=1}^{s_0}x_{u_i}\leq s_0x_{u_0}$.
Then $x_{u_{s_0+1}}\geq \frac{\rho-s_0}{\rho}x_{u_0}.$
Also, recall that $G'[S]\cong K_{1,e(S)}\cup (s-e(S)-1)K_1$. Then $y_{u_i}=y_{u_j}$ for $i,j \in [1,e(S)]$, and by considering the eigenvalue-eigenvector equation of  $A(G')$ with respect to $\rho'$, we obtain
$$\rho' y_{u_0}=\sum\limits_{i=1}^{e(S)}y_{u_i}+\sum_{i=0}^{t-1}y_{v_i}=e(S)\cdot y_{u_1}+\sum_{i=0}^{t-1}y_{v_i}~~\mbox{and}~~\rho' y_{u_1}=y_{u_0}+\sum_{i=0}^{t-1}y_{v_i}.$$
This implies that
$y_{u_1}=\frac{\rho'+1}{\rho'+e(S)}y_{u_0}\leq \frac{\rho+1}{\rho+e(S)}y_{u_0}$
because $\rho'\leq \rho$ and $e(S)>0$.
Thus,
$$y_{u_0}x_{u_{s_0+1}}-x_{u_0}y_{u_1}\geq x_{u_0}y_{u_0}\left(\frac{\rho-s_0}{\rho}-\frac{\rho+1}{\rho+e(S)}\right)
=x_{u_0}y_{u_0}\left(\frac{\rho(e(S)-(s_0+1))-s_0e(S)}{\rho(\rho+e(S))}\right)>0
$$
for sufficiently large $n$. Then $Y^t(\rho'-\rho)X\geq|E^\star|(y_{u_0}x_{u_{s_0+1}}-x_{u_0}y_{u_1})>0$, and so $\rho'>\rho$, a contradiction. Similarly, $e(T)\in\{t_0,t_0+1\}$.
This completes the proof of the claim.
\end{proof}

\begin{claim}\label{claim-k3}
If  $e(S)\neq 3$ (resp. $e(T)\neq 3$), then $e(S)=s_0$ and $G^\star[S]\cong K_{1,e(S)}\cup (s-e(S)-1)K_1$ (resp. $e(T)=t_0$ and $G^\star[T]\cong K_{1,e(T)}\cup (t-e(T)-1)K_1$), and if $e(S)=3$ (resp. $e(T)=3$), then $e(S)=s_0+1$ and $G^\star[S]\cong K_3\cup (s-3)K_1$ (resp. $e(T)=t_0+1$ and $G^\star[T]\cong K_3\cup (t-3)K_1$).
\end{claim}
\renewcommand\proofname{\bf Proof}
\begin{proof}
By symmetry, we only need to prove the claim for $S$. First assume that $e(S)\neq 3$.   By Claim \ref{claim-k2}, we have $e(S)\in \{s_0, s_0+1\}$.  If $e(S)=e(N_S[u_0])=s_0$, then the result follows. Thus we may assume that $e(S)=e(N_S[u_0])=s_0+1=d_S(u_0)+1=|N_s[u_0]|$. Since $e(S)\neq 3$, we claim that $s_0\geq 3$ and  $G^\star [N_S[u_0]]\cong K_{1,s_0}+u_1u_2$.  Let $G'=G-u_1u_2+u_0u_{s_0+1}$. It is easy to see that $G'[S]\cong K_{1,s_0+1}\cup (s-s_0-2)K_1$, and $G'$ is  $\Gamma_k$-free. Thus
$$
\begin{aligned}
Y^t(\rho'-\rho)X &=Y^t(A(G')-A(G^\star))X= x_{u_0}y_{u_{s_0+1}}+y_{u_0}x_{u_{s_0+1}}-x_{u_1}y_{u_2}-x_{u_2}y_{u_1}\\
&=x_{u_0}y_{u_1}+y_{u_0}x_{u_{s_0+1}}-2x_{u_1}y_{u_1},
\end{aligned}
$$
where the last equality follows from the fact that $x_{u_1}=x_{u_2}$ and $y_{u_i}=y_{u_j}$ for $i,j\in [1,s_0+1]$.
 Notice that $G^\star[S]\cong (K_{1,s_0}+u_1u_2)\cup (s-s_0-1)K_1$, $x_{u_1}=x_{u_2}$ and $x_{u_i}=x_{u_j}$ for $i,j \in [3,s_0].$ By considering the eigenvalue-eigenvector equation of  $A(G^\star)$ with respect to $\rho$, we obtain
$$
\begin{aligned}
\rho &x_{u_0}=x_{u_1}+x_{u_2}+\sum_{i=3}^{s_0}x_{u_i}+\sum_{i=0}^{t-1}x_{v_i}=2x_{u_1}+(s_0-2)x_{u_3}+\sum_{i=0}^{t-1}x_{v_i},\\
\rho &x_{u_1}=x_{u_0}+x_{u_2}+\sum_{i=0}^{t-1}x_{v_i}=x_{u_0}+x_{u_1}+\sum_{i=0}^{t-1}x_{v_i},~~
\rho x_{u_{3}}=x_{u_0}+\sum_{i=0}^{t-1}x_{v_i},~~\rho x_{u_{s_0+1}}=\sum_{i=0}^{t-1}x_{v_i},\\
\end{aligned}
$$
which leads to
$x_{u_1}=x_{u_2}=\frac{\rho(\rho+1)}{\rho^2+(s_0-1)\rho-(s_0-2)}x_{u_0}\leq \frac{\rho(\rho+1)}{\rho^2+2\rho-1}x_{u_0}$
and $x_{u_{s_0+1}}=\frac{\rho-1}{\rho}x_{u_1}-\frac{1}{\rho}x_{u_0}.$
Therefore,  $x_{u_0}y_{u_1}+y_{u_0}x_{u_{s_0+1}}-2x_{u_1}y_{u_1}$ is equal to
$$
\begin{aligned}
x_{u_0}y_{u_1}+y_{u_0}\left(\frac{\rho-1}{\rho}x_{u_1}-\frac{1}{\rho}x_{u_0}\right)-2x_{u_1}y_{u_1}
=x_{u_0}\left(y_{u_1}-\frac{y_{u_0}}{\rho}\right)-x_{u_1}\left(2y_{u_1}-\frac{\rho-1}{\rho}y_{u_0}\right).
\end{aligned}
$$
Recall that $G'[S]\cong K_{1,s_0+1}\cup (s-s_0-2)K_1$. According to  $\rho' y_{u_0}=(s_0+1)\cdot y_{u_1}+\sum_{i=0}^{t-1}y_{v_i}$,
$\rho' y_{u_1}=y_{u_0}+\sum_{i=0}^{t-1}y_{v_i}$  and the fact that $s_0\geq 3$, we get
\begin{equation*}
y_{u_1}=\frac{\rho'+1}{\rho'+s_0+1}y_{u_0}\leq \frac{\rho+1}{\rho+4}y_{u_0},
\end{equation*}
and therefore,
\begin{equation*}
2y_{u_1}-\frac{\rho-1}{\rho}y_{u_0}\geq \left(\frac{2(\rho'+1)}{\rho'+s_0+1}-1\right)y_{u_0}>0.
\end{equation*}
Notice that $x_{u_1}\leq \frac{\rho(\rho+1)}{\rho^2+2\rho-1}x_{u_0}.$
Then
$$
\begin{aligned}
&x_{u_0}\left(y_{u_1}-\frac{y_{u_0}}{\rho}\right)-x_{u_1}\left(2y_{u_1}-\frac{\rho-1}{\rho}y_{u_0}\right)
\geq
x_{u_0}\left(y_{u_1}-\frac{y_{u_0}}{\rho}\right)-\frac{\rho(\rho+1)}{\rho^2+2\rho-1}x_{u_0}\left(2y_{u_1}-\frac{\rho-1}{\rho}y_{u_0}\right)\\
&=\left(\frac{\rho^2-1}{\rho^2+2\rho-1}-\frac{1}{\rho}\right)x_{u_0}y_{u_0}-\frac{\rho^2+1}{\rho^2+2\rho-1}x_{u_0}y_{u_1}\\
&\geq \left(\frac{\rho^2-1}{\rho^2+2\rho-1}-\frac{1}{\rho}-\frac{(\rho^2+1)(\rho+1)}{(\rho^2+2\rho-1)(\rho+4)}\right)x_{u_0}y_{u_0}>0.
\end{aligned}
$$
Thus $Y^t(\rho'-\rho)X =x_{u_0}y_{u_1}+y_{u_0}x_{u_{s_0+1}}-2x_{u_1}y_{u_1}=x_{u_0}\left(y_{u_1}-\frac{y_{u_0}}{\rho}\right)-x_{u_1}\left(2y_{u_1}-\frac{\rho-1}{\rho}y_{u_0}\right) >0,$ implying that $ \rho'>\rho$, a contradiction.
Therefore, if $e(S)\neq 3$, then $e(S)=s_0$ and $G^\star[S]\cong K_{1,e(S)}\cup (s-e(S)-1)K_1$.

Now suppose that $e(S)=3$. If $e(S)=e(N_S[u_0])=s_0+1$, then $s_0=2$ and $G^\star[S]\cong K_3\cup (s-3)K_1$, as required. So the remaining case is $e(S)=e(N_S[u_0])=s_0$. In this situation, we have $s_0=3$ and $G^\star[S]\cong K_{1,3}\cup (s-4)K_1$. Let $G'=G^\star-u_0u_3+u_1u_2$. Then it is clear that $G'[S]\cong K_3\cup (s-3)K_1$, and
$$
\begin{aligned}
Y^t(\rho'-\rho)X &= Y^t(A(G')-A(G^\star))X = x_{u_1}y_{u_2}+x_{u_2}y_{u_1}-x_{u_0}y_{u_3}-x_{u_3}y_{u_0}
  =x_{u_1}y_{u_0}-x_{u_0}y_{u_3},
  \end{aligned}
$$
where the last equality follows from the fact that $x_{u_1}=x_{u_2}=x_{u_3}$ and $y_{u_0}=y_{u_1}=y_{u_2}$.
As above, by $\rho x_{u_0}=3x_{u_1}+\sum\limits_{i=0}^{t-1}x_{v_i}$, $\rho x_{u_1}=x_{u_0}+\sum\limits_{i=0}^{t-1}x_{v_i}$,
$\rho' y_{u_0}=2y_{u_0}+\sum\limits_{i=0}^{t-1}y_{v_i}$ and $\rho'y_{u_3}=\sum\limits_{i=0}^{t-1}y_{v_i}$,
we obtain
\begin{equation*}
x_{u_1}=\frac{\rho+1}{\rho+3}x_{u_0}~~\mbox{and}~~
y_{u_3}=\frac{\rho'-2}{\rho'}y_{u_0}\leq \frac{\rho-2}{\rho}y_{u_0}.
 \end{equation*}
 Hence,
 \begin{equation*}
  x_{u_1}y_{u_0}-x_{u_0}y_{u_3}\geq x_{u_0}y_{u_0}\left(\frac{\rho+1}{\rho+3}-\frac{\rho-2}{\rho}\right)
   = \frac{6}{\rho(\rho+3)}x_{u_0}y_{u_0} >0.
   \end{equation*}
Then $Y^t(\rho'-\rho)X =x_{u_1}y_{u_0}-x_{u_0}y_{u_3}>0,$
and so $\rho'>\rho$, contrary to the maximality of $\rho(G^\star)$. Hence, if $e(S)=3$, then $e(S)=s_0+1$ and $G^\star[S]\cong K_3\cup (s-3)K_1$.
\end{proof}

According to Claim \ref{claim-k3}, we only need to consider the following four cases.
For simplicity, we denote by $s^\star=e(S)$ and $t^\star=e(T)$.

\noindent{\bf{$\underline{\mbox{Case 1.}}$}}  $G^\star[S] \cong K_{3}\cup (s-3)K_1$ and $G^\star[T] \cong K_{3}\cup (t-3)K_1$.

In this situation, we have $x_{u_0}=x_{u_1}=x_{u_2}$, $x_{u_i}=x_{u_j}$ for $i,j \in [3, s-1]$ and $x_{v_0}=x_{v_1}=x_{v_2},$ $x_{v_i}=x_{v_j}$ for $i,j \in [3, t-1]$.
Combining $\rho x_{u_0}=x_{u_1}+x_{u_2}+\sum\limits_{i=0}^{t-1}x_{v_i}=2 x_{u_0}+\rho x_{u_3}$ and $\rho x_{v_0}=x_{v_1}+x_{v_2}+\sum\limits_{i=0}^{s-1}x_{u_i}=2 x_{v_0}+\rho x_{v_3}$ yields that
\begin{equation}\label{001}
 x_{v_3}=\frac{\rho-2}{\rho}x_{v_0}, ~x_{u_3} = \frac{\rho-2}{\rho}x_{u_0}~~\mbox{and}~~
  (\rho-2)(x_{v_0}-x_{u_0})=\sum\limits_{i=0}^{s-1}x_{u_i}-\sum\limits_{i=0}^{t-1}x_{v_i}.
\end{equation}
Furthermore, we assert that $x_{v_0}\leq x_{u_0}$. Suppose to the contrary that $x_{v_0}> x_{u_0}$, then  we obtain $\sum\limits_{i=0}^{s-1}x_{u_i}>\sum\limits_{i=0}^{t-1}x_{v_i}$.
On the other hand,
\begin{equation*}
\begin{aligned}
\sum_{i=0}^{t-1}x_{v_i} &= \sum_{i=0}^2 x_{v_i}+\sum_{i=3}^{t-1}x_{v_i}=3x_{v_0}+(t-3)x_{v_3}= 3x_{v_0}+(t-3)\frac{\rho-2}{\rho}x_{v_0}\\
&>3x_{u_0}+(s-3)\frac{\rho-2}{\rho}x_{u_0}= \sum_{i=0}^{s-1}x_{u_i},
\end{aligned}
\end{equation*}
which is a  contradiction.

Let $E^\star=\{v_0v_1,v_0v_2,v_1v_2,\}$, $E'=\{u_0u_3,u_0u_4,u_3u_4\}$ and $G'=G^\star-E^\star+E'$. We see that $G'[S]\cong (K_1\nabla 2K_2)\cup (s-5)K_1$ and $G'$ is $\Gamma_k$-free, where $K_1\nabla 2K_2$ is a graph obtained from the disjoint union $K_1 \cup 2K_2$ by adding all edges between $K_1$ and $2K_2.$
Then $y_{u_1}=y_{u_2}=y_{u_3}=y_{u_4}$, $y_{u_i}=y_{u_j}$ for $i,j \in [5, s-1]$ and $y_{v_i}=y_{v_j}$ for $i,j \in [0, t-1]$. By considering the eigenvalue-eigenvector equation of  $A(G')$ with respect to $\rho'$, we obtain
$$
\begin{aligned}
     \rho' y_{u_0}=4y_{u_1}+\sum\limits_{i=0}^{t-1}y_{v_i},~
     \rho' y_{u_1}=y_{u_0}+y_{u_1}+\sum\limits_{i=0}^{t-1}y_{v_i}, ~
      \rho' y_{u_5}=\sum\limits_{i=0}^{t-1}y_{v_i} ,~~
      \rho' y_{v_0}=y_{u_0}+4y_{u_1}+(s-5)y_{u_5},
\end{aligned}
$$
which gives that
\begin{equation}\label{1.}
 y_{u_5}=\left(1-\frac{4(\rho'+1)}{\rho'(\rho'+3)}\right)y_{u_0}~\mbox{and}~y_{u_1}=\frac{\rho'+1}{\rho'+3}y_{u_0}.
\end{equation}
Hence,
 $y_{v_0} =\frac{y_{u_0}}{\rho'}\left(s-4+\frac{4(\rho'+1)(\rho'-s+5)}{\rho'(\rho'+3)}\right)
=\frac{y_{u_0}}{\rho'}\left(s-4+\frac{4(\rho'+3)^2-4s(\rho'+3)+8(s-2)}{\rho'(\rho'+3)}\right)$.  Combining this with $\rho'> s=\lfloor\frac{n}{2}\rfloor$ (since $G'$ contains  $K_{\lfloor\frac{n}{2}\rfloor,\lceil\frac{n}{2}\rceil}$ as a proper subgraph), we deduce that
\begin{equation}\label{1.1}
y_{v_0} \leq \frac{y_{u_0}}{\rho'}\left(\rho'-4+\frac{4(\rho'+3)^2-4\rho'(\rho'+3)+8(\rho'-2)}{\rho'(\rho'+3)}\right)
           =\frac{y_{u_0}}{\rho'}(\rho'-4)+\emph{O}\left(\frac{1}{\rho'^2}\right).
\end{equation}
Recall that $x_{v_0}\leq x_{u_0}$, and  $X$ and $Y$ are the Perron vector of $G^\star$ and $G'$, respectively. Combining  (\ref{001}), (\ref{1.}) and (\ref{1.1}), we have
\begin{equation*}
\begin{aligned}
Y^t(\rho'-\rho)X&=Y^t(A(G')-A(G))X= \sum\limits_{ij\in E'}(x_{i}y_{j}+y_{i}x_{j})-\sum\limits_{ij\in E^\star}(x_{i}y_{j}+y_{i}x_{j})\\
&=2((y_{u_1}x_{u_0}+(y_{u_0}+y_{u_1})x_{u_3})-3x_{v_0}y_{v_0})\geq 2((y_{u_1}x_{u_0}+(y_{u_0}+y_{u_1})x_{u_3})-3x_{u_0}y_{v_0})\\
&=2x_{u_0}y_{u_0}\left(\frac{\rho'+1}{\rho'+3}+(1+\frac{\rho'+1}{\rho'+3})\frac{\rho-2}{\rho}\right)-6x_{u_0}y_{v_0}\\
&\geq2x_{u_0}y_{u_0}\left(\frac{\rho'+1}{\rho'+3}+(1+\frac{\rho'+1}{\rho'+3})\frac{\rho'-2}{\rho'}\right)-6x_{u_0}y_{v_0}\\
&\geq 2x_{u_0}y_{u_0}\left(3-\frac{4}{\rho'}-\frac{4}{\rho'+3}-         \frac{3}{\rho'}(\rho'-4)+\emph{O}\left(\frac{1}{\rho'^2}\right)\right)\\
&=2x_{u_0}y_{u_0}\left(\frac{12}{\rho'}-\frac{4}{\rho'}-\frac{4}{\rho'+3}+\emph{O}\left(\frac{1}{\rho'^2}\right)\right)>0.
\end{aligned}
\end{equation*}
It follows that $\rho'>\rho$, a contradiction.

\noindent{\bf{$\underline{\mbox{Case 2.}}$}}  $G^\star[S] \cong K_{1,s^\star}\cup (s-s^\star-1)K_1$ and $G^\star[T] \cong K_{1,t^\star}\cup (t-t^\star-1)K_1$, where $s^\star\neq 3$ and $t^\star\neq 3$.

Note that $s_0=e(N_S[u_0])$, $t_0=e(N_T[v_0])$ and  $s_0+t_0=k-1$.  If $s_0=k-1,$ then $G^\star \cong K_{\lceil\frac{n}{2}\rceil,\lfloor\frac{n}{2}\rfloor}\diamond K_{1,k-1},$ as desired. So we may assume that  $s_0\leq k-2.$
Observe that $x_{u_i}=x_{u_j}$ for $i,j\in [1,s_0]$ and  $x_{v_i}=x_{v_j}$ for $i,j\in [1,t_0]$.
Let $E^\star=\{v_0v_i: 1\leq i \leq t_0\}$, $E'=\{u_0u_i: s_0+1\leq i\leq k-1\}$, and $G'=G^\star-E^\star+E'.$
It is easy to see that $G'[S]\cong  K_{1,k-1}\cup (s-k)K_1$ and $G'$ is $\Gamma_k$-free.
By symmetry, $y_{u_i}=y_{u_j}$ for $i,j\in[1,k-1]$. Then
$$
\begin{aligned}
\rho y_{u_0}=\sum_{i=1}^{k-1}y_{u_i}+\sum_{i=0}^{t-1}y_{v_i}=(k-1)y_{u_1}+\sum_{i=0}^{t-1}y_{v_i},~
\rho y_{u_1}=y_{u_0}+\sum_{i=0}^{t-1}y_{v_i},~
\rho y_{u_{k}}=\sum_{i=0}^{t-1}y_{v_i},
\end{aligned}
$$
which leads to
$$ y_{u_1}=\frac{\rho'+1}{\rho'+k-1}y_{u_0}, ~y_{u_k}=\frac{\rho'y_{u_1}-y_{u_0}}{\rho'},$$
and
$$ \sum\limits_{i=0}^{s-1}y_{u_i}=\frac{\rho'^2s+2\rho' (k-1)-(k-1)s+k(k-1)}{\rho'(\rho'+k-1)}y_{u_0}\leq \frac{\rho'^3+\rho'(k-1)+k(k-1)}{\rho'(\rho'+k-1)}y_{u_0}.$$
Notice that $\rho'y_{v_0}=\sum\limits_{i=0}^{s-1}y_{u_i}.$ Then
\begin{equation}\label{case2-2}
  y_{v_0}\leq\frac{\rho'^3+\rho'(k-1)+k(k-1)}{\rho'^2(\rho'+k-1)}y_{u_0}.
\end{equation}
Similarly, we have
$x_{u_1}=\frac{\rho+1}{\rho+s_0}x_{u_0}$,
$x_{u_{s_0+1}}=\frac{\rho^2-s_0}{\rho(\rho+s_0)}x_{u_0}$. Since $s=\lfloor\frac{n}{2}\rfloor< \rho$, it follows that
$$\sum\limits_{i=0}^{s-1}x_{u_i}=\frac{\rho^2s+2\rho s_0-s_0s+s_0^2+s_0}{\rho(\rho+s_0)}x_{u_0}\leq \frac{\rho^3+\rho s_0+s_0^2+s_0}{\rho(\rho+s_0)}x_{u_0}.$$
Hence, $x_{v_1}=\frac{\rho+1}{\rho+t_0}x_{v_0}$.
According to  $\rho x_{v_0}=\sum\limits_{i=1}^{t_0}x_{v_i}+\sum\limits_{i=0}^{s-1}x_{u_i}=t_0x_{v_1}+\sum\limits_{i=0}^{s-1}x_{u_i},$
we obtain
\begin{equation*}
  x_{v_0}\leq\frac{(\rho+k-1-s_0)(\rho^3+\rho s_0+s_0^2+s_0)}{\rho(\rho+s_0)(\rho^2-(k-1-s_0))}x_{u_0}.
\end{equation*}
It follows that
\begin{equation}\label{case2-2.}
\begin{aligned}
 &\frac{2\rho+k-s_0}{\rho+k-1-s_0}x_{v_0}y_{v_0}\\
 &\leq  \frac{2\rho+k-s_0}{\rho+k-1-s_0}\cdot \frac{(\rho+k-1-s_0)(\rho^3+\rho s_0+s_0^2+s_0)}{\rho(\rho+s_0)(\rho^2-(k-1-s_0))}\cdot\frac{\rho'^3+\rho'(k-1)+k(k-1)}  {\rho'^2(\rho'+k-1)}x_{u_0}y_{u_0} \\
      &= \frac{(2\rho+k-s_0)(\rho^3+\rho s_0+s_0(s_0+1))(\rho'^3+\rho'(k-1)+k(k-1))}{\rho \rho'^2(\rho+s_0)(\rho^2-(k-1-s_0))(\rho'+k-1)}x_{u_0}y_{u_0}\\
      &=\left(\frac{2\rho^4\rho'^3+\rho^3\rho'^3(k-s_0)}{\rho \rho'^2(\rho+s_0)(\rho^2-(k-1-s_0))(\rho'+k-1)}
      +\emph{O}\left(\frac{1}{\rho^2}\right)\right)x_{u_0}y_{u_0}.
\end{aligned}
\end{equation}
Note that $X$ and $Y$ are the Perron vectors of $G^\star$ and $G'$, respectively. Then
\begin{equation}\label{case2-3}
  \begin{aligned}
  Y^t(\rho'-\rho)X&=Y^t(A(G')-A(G^\star))X=\sum\limits_{u_0u_i\in E'}(x_{u_0}y_{u_i}+y_{u_0}x_{u_i})-\sum\limits_{v_0v_i\in E^\star}(x_{v_0}y_{v_i}+y_{v_0}x_{v_i})\\
                  &=(k-1-s_0)(x_{u_0}y_{u_1}+y_{u_0}x_{u_{s_0+1}}-x_{v_0}y_{v_0}-x_{v_1}y_{v_0})\\
                  &=(k-1-s_0)\left(\left(\frac{\rho'+1}{\rho'+k-1}+\frac{\rho^2-s_0}{\rho(\rho+s_0)}\right)x_{u_0}y_{u_0}
                  -\frac{2\rho+k-s_0}{\rho+k-1-s_0}x_{v_0}y_{v_0}\right).
  \end{aligned}
\end{equation}
Recall that $s_0\leq k-2$. We shall prove $(\ref{case2-3})>0$ by showing $$\left(\frac{\rho'+1}{\rho'+k-1}+\frac{\rho^2-s_0}{\rho(\rho+s_0)}\right)x_{u_0}y_{u_0}
>\frac{2\rho+k-s_0}{\rho+k-1-s_0}x_{v_0}y_{v_0},$$ which leads to $\rho'>\rho$, 
and we derive a contradiction.
According to (\ref{case2-2.}), it suffices to show
$$\frac{\rho'+1}{\rho'+k-1}+\frac{\rho^2-s_0}{\rho(\rho+s_0)}
                  >\frac{2\rho^4\rho'^3+\rho^3\rho'^3(k-s_0)}{\rho \rho'^2(\rho+s_0)(\rho^2-(k-1-s_0))(\rho'+k-1)}
      +\emph{O}(\frac{1}{\rho^2}).$$
      The last inequality holds by $\frac{\rho'+1}{\rho'+k-1}+\frac{\rho^2-s_0}{\rho(\rho+s_0)}=\frac{2\rho^4(\rho'^3+k\rho^4\rho'^2+s_0\rho^3\rho'^3}{\rho \rho'^2(\rho+s_0)(\rho^2-(k-1-s_0))(\rho'+k-1)}
      +\emph{O}\left(\frac{1}{\rho^2}\right),$ as required.

\noindent{\bf{$\underline{\mbox{Case 3.}}$}}  $G^\star[S] \cong K_{1,s^\star}\cup (s-s^\star-1)K_1$ and $G^\star[T] \cong K_3\cup (t-3)K_1$, where $s^\star\neq 3$.

Clearly, $k\geq 4$ and $s_0=k-4$. First we may assume that $k\geq 5$. Then $x_{u_i}=x_{u_j}$ for $i,j \in [1,k-4]$ and $x_{v_0}=x_{v_1}=x_{v_2}$. It follows that
\begin{equation}\label{iii-0}
  x_{u_1}=\frac{\rho+1}{\rho+k-4}x_{u_0},
\end{equation}
and
\begin{equation}\label{iii-1}
 x_{u_{k-3}}=\frac{\rho^2-k+4}{\rho(\rho+k-4)}x_{u_0}=\left(\frac{\rho}{\rho+k-4}+\emph{O}\left(\frac{1}{\rho^2}\right)\right)x_{u_0}
  =\left(1-\frac{k-4}{\rho+k-4}+\emph{O}\left(\frac{1}{\rho^2}\right)\right)x_{u_0}.
\end{equation}
Combining with $s=\lfloor\frac{n}{2}\rfloor<\rho$, we have
\begin{equation*}
    \begin{aligned}
    \sum\limits_{i=0}^{s-1}x_{u_i} &=x_{u_0}+(k-4)x_{u_1}+(s-k+3)x_{u_{k-3}}
         = \frac{\rho^2s+2\rho(k-4)-s(k-4)+(k-4)(k-3)}{\rho(\rho+k-4)}x_{u_0} \\
        &\leq \frac{\rho^3+2\rho(k-4)-\rho(k-4)+(k-4)(k-3)}{\rho(\rho+k-4)}x_{u_0}
        =\left(\frac{\rho^2+k-4}{\rho+k-4}+\emph{O}\left(\frac{1}{\rho^2}\right)\right)x_{u_0}.
      \end{aligned}
\end{equation*}
Notice that $\rho x_{v_0}=x_{v_1}+x_{v_2}+\sum\limits_{i=0}^{s-1}x_{u_i}=2x_{v_0}+\sum\limits_{i=0}^{s-1}x_{u_i}.$
Then we obtain 
\begin{equation}\label{iii-3}
  x_{v_0}\leq \frac{1}{\rho-2}\left(\frac{\rho^2+k-4}{\rho+k-4}+\emph{O}\left(\frac{1}{\rho^2}\right)\right)x_{u_0}
         =\left(\frac{\rho^2+k-4}{(\rho-2)(\rho+k-4)}+\emph{O}\left(\frac{1}{\rho^3}\right)\right)x_{u_0}.
\end{equation}
Let $E^\star=\{v_0v_1, v_0v_2, v_1v_2\}$, $E'=\{u_0u_{k-3}, u_0u_{k-2}, u_0u_{k-1}\},$ and $G'=G^\star-E^\star+E'$.
Then $G'\cong K_{\lceil\frac{n}{2}\rceil,\lfloor\frac{n}{2}\rfloor}\diamond K_{1,k-1}$ and $G'$ is $\Gamma_k$-free. Note that
\begin{equation}\label{iii-0'}
  y_{u_1}=\frac{\rho'+1}{\rho'+k-1}y_{u_0}=(1-\frac{k-2}{\rho'+k-1})y_{u_0}, ~y_{u_k}=\frac{\rho'y_{u_1}-y_{u_0}}{\rho'},
\end{equation}
and
\begin{equation*}
\begin{aligned}
\sum\limits_{i=0}^{s-1}y_{u_i} &=\frac{\rho'^2s+2\rho'(k-1)-s(k-1)+k(k-1)}{\rho'(\rho'+k-1)}y_{u_0}\\
&\leq\frac{\rho'^3+2\rho'(k-1)-\rho'(k-1)+k(k-1)}{\rho'(\rho'+k-1)}y_{u_0}=\left(\frac{\rho'^2+k-1}{\rho'+k-1}+\emph{O}\left(\frac{1}{\rho'^2}\right)\right)y_{u_0},
\end{aligned}
\end{equation*}
where the inequality by $s=\lfloor\frac{n}{2}\rfloor<\rho'$.
Then by using $\rho' y_{v_0}=\sum\limits_{i=0}^{s-1}y_{u_i}$, we have
\begin{equation}\label{iii-3'}
\begin{aligned}
y_{v_0}\leq \frac{1}{\rho'}\left(\frac{\rho'^2+k-1}{\rho'+k-1}+\emph{O}\left(\frac{1}{\rho'^2}\right)\right)y_{u_0}
          =\left(\frac{\rho'^2+k-1}{\rho'(\rho'+k-1)}+\emph{O}\left(\frac{1}{\rho'^3}\right)\right)y_{u_0}.
\end{aligned}
\end{equation}
Thus, combining this with (\ref{iii-3}), we deduce
\begin{equation}\label{iii-4}
\begin{aligned}
\vspace{3mm}
   x_{v_0}y_{v_0} & \leq \left(\frac{\rho^2+k-4}{(\rho-2)(\rho+k-4)}+\emph{O}\left(\frac{1}{\rho^3}\right)\right)
  \left(\frac{\rho'^2+k-1}{\rho'(\rho'+k-1)}+\emph{O}\left(\frac{1}{\rho'^3}\right)\right)x_{u_0}y_{u_0}\\
                  &=\left(\frac{(\rho^2+k-4)(\rho'^2+k-1)}{\rho'(\rho'+k-1)(\rho-2)(\rho+k-4)}
                  +\emph{O}\left(\frac{1}{\rho^3}\right)\right)x_{u_0}y_{u_0}.
\end{aligned}
\end{equation}
By symmetry, we have
\begin{equation}\label{iii-5}
\begin{aligned}
Y^t(\rho'-\rho)X&=Y^t(A'-A)X=\sum\limits_{ij\in E'}(x_iy_j+x_jy_i)-\sum\limits_{ij\in E^\star}(x_iy_j+x_jy_i)\\
   &=3x_{u_0}y_{u_1}+3y_{u_0}x_{u_{k-3}}-6x_{v_0}y_{v_0}.
\end{aligned}
\end{equation}
Now, we shall show $(\ref{iii-5})>0$. According to (\ref{iii-1}) and (\ref{iii-0'}), we change (\ref{iii-5}) to
$3x_{u_0}y_{u_0}(\frac{\rho'+1}{\rho'+k-1}+\frac{\rho^2-k+4}{\rho(\rho+k-4)})
   -6x_{v_0}y_{v_0}.$
Thus, combining with (\ref{iii-4}), it suffices to prove
$$\frac{\rho'+1}{\rho'+k-1}+\frac{\rho^2-k+4}{\rho(\rho+k-4)}
>2\left(\frac{(\rho^2+k-4)(\rho'^2+k-1)}{\rho'(\rho'+k-1)(\rho-2)(\rho+k-4)}+\emph{O}\left(\frac{1}{\rho^3}\right)\right).$$
Let us multiply both sides of this inequality by $\rho'\rho(\rho'+k-1)(\rho-2)(\rho+k-4).$
Then it suffices to show $$\rho\rho'(\rho'+1)(\rho-2)(\rho+k-4)+\rho'(\rho'+k-1)(\rho-2)(\rho^2-k+4)>2\rho(\rho^2+k-4)(\rho'^2+k-1)+\emph{O}(\rho^2).$$
By calculation, we only need to prove that $\rho^2\rho'^2(k-8)+\rho^3\rho'k>0$ since $\rho= \emph{O}(n), \rho'= \emph{O}(n).$
Note that $\rho>\rho'$. Then $$\rho^2\rho'^2(k-8)+\rho^3\rho'k>\rho^2\rho'^2(k-8)+\rho^2\rho'^2k=\rho^2\rho'^2(2k-8).$$ Thus, for $k\geq 5$, we have $\rho'>\rho$ by (\ref{iii-5}), a contradiction.

For $k=4,$ we see that $G^\star[N_T[v_0]] \cong K_{3}$, and $ G^\star[N_S[u_0]]$ is an empty graph.
If $s=t=\frac{n}{2}$, then $G^\star\cong K_{\frac{n}{2},\frac{n}{2}}\diamond K_3$, as desired.
If $s=t-1=\frac{n-1}{2}$, then let $E^\star=\{v_0v_1, v_0v_2, v_1v_2\}$, $E'=\{u_0u_1, u_0u_2, u_1u_2 \}$ and $G'=G^\star-E^\star+E'.$
Then $G' \cong K_{\lceil\frac{n}{2}\rceil,\lfloor\frac{n}{2}\rfloor}\diamond K_3$ and $G'$ is $\Gamma_4$-free.
Notice that the quotient matrix of $A(G^\star)$ with respect to the partition $\Pi^\star: V(G^\star)=\{u_0,\cdots, u_{s-1}\}\cup \{v_0, v_1, v_2\}\cup \{v_3,\cdots, v_{t-1}\}$ is given by
$$
B_1=\left(
  \begin{array}{ccc}
    0 & 3 & t-3 \\
    s & 2 & 0 \\
    s & 0 & 0 \\
  \end{array}
\right).
$$
Furthermore, the characteristic polynomial of  $B_1$ is
$\varphi(B_1,x)=x^3-2x^2-stx+2s(t-3),$ and its  largest root coincides with $\rho$ by Lemma \ref{quotient} since the partition is equitable.
By symmetry, the characteristic polynomial  $\varphi(B_2,x)$ of the quotient matrix  $B_2$ of $A(G')$ with respect to the partition $\Pi': V(G') =\{v_0,\cdots, v_{t-1}\}\cup \{u_0, u_1, u_2\}\cup \{u_3,\cdots, u_{s-1}\}$ can be obtained from $\varphi(B_1,x)$ by switching $s$ and $t$. 
The partition is also equitable, then by Lemma \ref{quotient}, $\rho'=\rho(B_2)$.
Notice that $\rho>s=\frac{n-1}{2}$ and $s+t=n.$ By a simple calculation, we have $$\varphi(B_2,\rho) = \varphi(B_2,\rho)-\varphi(B_1,\rho) = -4\rho^2-\left(\frac{n^2}{2}-\frac{1}{2}\right)\rho+n^2-6n-1. $$
Thus $\varphi(B_2,\rho)<\varphi(B_2,\frac{n-1}{2})=-\frac{n^3}{4}+\frac{n^2}{4}-\frac{15n}{4}-\frac{9}{4}<0$, as required. Therefore, $\rho<\rho'$, which leads to a contradiction.

\noindent{\bf{$\underline{\mbox{Case 4.}}$}}  $G^\star[S] \cong K_3\cup (s-3)K_1$ and $G^\star[T] \cong K_{1,t^\star}\cup (t-t^\star-1)K_1$, where $t^\star\neq 3$.

Then $x_{u_0}=x_{u_1}=x_{u_2}$ and $x_{v_i}=x_{v_j}$ for $i,j \in[1,k-4]$.
It follows that
\begin{equation}\label{iii-00}
\begin{aligned}
    x_{v_1}=\frac{\rho+1}{\rho+k-4}x_{v_0},~x_{u_{3}}=\frac{\rho-2}{\rho}x_{u_0},
\end{aligned}
\end{equation}
and
\begin{equation*}
\begin{aligned}
  \sum\limits_{i=0}^{s-1} x_{u_i} =3x_{u_0}+(s-3)x_{u_3}= \frac{\rho s-2s+6}{\rho}x_{u_0} \leq \frac{\rho^2+2\rho+6}{\rho}x_{u_0}=(\rho-2+\frac{6}{\rho})x_{u_0},
\end{aligned}
\end{equation*}
where the inequality above holds as $s=\lfloor\frac{n}{2}\rfloor<\rho$. Notice that $$\rho x_{v_0}=\sum\limits_{i=1}^{k-4}x_{v_i}+\sum\limits_{i=0}^{s-1} x_{u_i} =(k-4)x_{v_1}+\sum\limits_{i=0}^{s-1} x_{u_i}.$$ Then
\begin{equation}\label{iii-33}
\begin{aligned}
   x_{v_0} = \left(\frac{(k-4)(\rho+1)}{\rho(\rho+k-4)}+\frac{\rho-2}{\rho}+\frac{6}{\rho^2}\right)x_{u_0}
   \leq \left(1+\frac{k-4}{\rho+k-4}-\frac{2}{\rho}+\emph{O}\left(\frac{1}{\rho^2}\right)\right)x_{u_0}.
\end{aligned}
\end{equation}
Let $E^\star=\{v_0v_1, v_0v_2,..., v_0v_{k-4},u_1u_2\}$, $E'=\{u_0u_{3}, u_0u_{4},..., u_0u_{k-1}\}$ and $G'=G^\star-E^\star+E'.$
Then $G'\cong K_{\lceil\frac{n}{2}\rceil,\lfloor\frac{n}{2}\rfloor}\diamond K_{1,k-1}$ and $G'$ is $\Gamma_k$-free.
By  (\ref{iii-3'}) and (\ref{iii-33}), we have
\begin{equation*}
\begin{aligned}
\vspace{3mm}
   x_{v_0}y_{v_0} & \leq \left(1-\frac{2}{\rho}+\frac{k-4}{\rho+k-4}+\emph{O}\left(\frac{1}{\rho^2}\right)\right)
  \left(\frac{\rho'}{\rho'+k-1}+\emph{O}\left(\frac{1}{\rho^2}\right)\right)x_{u_0}y_{u_0} \\
   \vspace{3mm}
   &=\left(1+\frac{k-4}{\rho+k-4}-\frac{2}{\rho}-\frac{k-1}{\rho'+k-1} +\emph{O}\left(\frac{1}{\rho^2}\right)\right)x_{u_0}y_{u_0}.
\end{aligned}
\end{equation*}
Combining with (\ref{iii-0'}) and (\ref{iii-00}), we have
\begin{equation*}
\begin{aligned}
   &Y^t(\rho'-\rho)X=Y^t(A'-A)X =\sum\limits_{ij\in E'}(x_iy_j+x_jy_i)-\sum\limits_{ij\in E^\star}(x_iy_j+x_jy_i)\\
 \vspace{3mm}
   &=(k-3)(x_{u_0}y_{u_1}+x_{u_1}y_{u_3})-(k-4)(x_{v_0}y_{v_1}+y_{v_0}x_{v_1})-2x_{u_0}y_{u_1}\\
 \vspace{3mm}
   &=((k-5)(1-\frac{k-2}{\rho'+k-1})+(k-3)(1-\frac{2}{\rho}))x_{u_0}y_{u_0}-(k-4)(2-\frac{k-5}{\rho+k-4})x_{v_0}y_{v_0}\\
  \vspace{3mm}
   &\geq\left(2k-8-\frac{(k-5)(k-2)}{\rho'+k-1}-\frac{2(k-3)}{\rho}\right)x_{u_0}y_{u_0}\\
   \vspace{3mm}
   &~~~~~-\left(2k-8-\frac{(k-5)(k-4)}{\rho+k-4}\right)
  \left(1+\frac{k-4}{\rho+k-4}-\frac{2}{\rho}-\frac{k-1}{\rho'+k-1}+\emph{O}\left(\frac{1}{\rho^2}\right)\right)x_{u_0}y_{u_0}\\
   &=\left(\frac{k^2-k-2}{\rho'+k-1}+\frac{2(k-5)}{\rho}-\frac{(k-4)(k-3)}{\rho+k-4} +\emph{O}\left(\frac{1}{\rho^2}\right)\right)x_{u_0}y_{u_0}>0.
\end{aligned}
\end{equation*}
It follows that $\rho'>\rho$, which is impossible. This completes the proof of Case 4.

Considering Case 1-4, we complete the proof of Theorem \ref{thm-2}. $\hfill\blacksquare$

\section{Conclusion remark}


Bollob\'{a}s \cite{B78} asked: what is the maximum size of  a graph $G\in \overline{\Omega'_k}$ with order $n$? In this paper, we pose a spectral analogue problem as follows.

\begin{prob}
What is the maximum spectral radius of a graph $G\in \overline{\Omega'_k}$ with $k\geq 3$ and order $n$?
\end{prob}

Ma and Yang \cite{MY} proved that $f(n)<n+\sqrt{n}+o(n)$ for any $n$-vertex 2-connected graph.
Theorem \ref{thm1.2} shows that the graph with the maximum spectral radius among all graphs with no two cycles of the same length has a cut vertex. So it is natural to ask the following problem.

\begin{prob}
What is the maximum spectral radius among all 2-connected $n$-vertex graphs with no two cycles of the same length?
\end{prob}

\end{document}